# Finding a maximal element of a convex set through its characteristic cone: An application to finding a strictly complementary solution

Mahmood Mehdiloozad[*]

*Department of Mathematics, College of Sciences, Shiraz University, Shiraz, Iran*

Kaoru Tone

*National Graduate Institute for Policy Studies, Tokyo, Japan*

Rahim Askarpour

*Department of Mathematics, College of Sciences, Shiraz University, Shiraz, Iran*

Mohammad Bagher Ahmadi

*Department of Mathematics, College of Sciences, Shiraz University, Shiraz, Iran*

[*] **Corresponding author**: M. Mehdiloozad

Ph.D. student

Department of Mathematics

College of Sciences

Shiraz University

Golestan Street | Adabiat Crossroad | Shiraz 71454 | Iran

E-mail: m.mehdiloozad@gmail.com

Tel.: + 98.9127431689

# Finding a maximal element of a convex set through its characteristic cone: An application to finding a strictly complementary solution

**Abstract** In order to express a polyhedron as the (Minkowski) sum of a polytope and a polyhedral cone, Motzkin (1936) made a transition from the polyhedron to a polyhedral cone. Based on his excellent idea, we represent a set by a *characteristic cone*. By using this representation, we then reach four main results: (i) expressing a closed convex set containing no line as the direct sum of the convex hull of its extreme points and conical hull of its extreme directions, (ii) establishing a convex programming (CP) based framework for determining a *maximal element*—an element with the maximum number of positive components—of a convex set, (iii) developing a linear programming problem for finding a relative interior point of a polyhedron, and (iv) proposing two procedures for the identification of a strictly complementary solution in linear programming.

## 1 Introduction

One of the most important theoretical results in linear programming (LP) is the *Goldman–Tucker theorem* (Goldman and Tucker, 1956) that establishes a significant relation between a pair of primal and dual LP problems. This relation, referred to as *strictly complementary slackness* *condition (SCSC)*[1], implies the existence of a pair of primal–dual optimal solutions for which the sum of each pair of complementary variables is positive. It should be noted here that this positivity is not guaranteed by the ordinary *complementary slackness condition (CSC)*. A pair of solutions satisfying the SCSC, for brevity, is called a *strictly complementary solution*. Though the existence of such a solution was originally proved by Goldman and Tucker (Goldman and Tucker, 1956) using the self-dual theorem of Tucker (1956), an elegant and constructive proof was provided later by Balinski and Tucker (1969). Another proof based on the interior point methods can also be found in Adler (1992) and Jansen (1994).

Under the existence of a strictly complementary solution, the issue of identifying this solution is of great importance for three reasons. First, this solution induces unique optimal partitions for the set of the indices of non-negative variables. These partitions are directly related to the behavior of the optimal objective value and are hence important in sensitivity analysis (Greenberg, 1994). Second, a strictly complementary solution determines the dimensions and the degeneracy degrees of the optimal faces of the primal and dual LP problems (Tijssen and Sierksma, 1998). Third, some successes are recently achieved in the appropriate use of a strictly complementary solution in the theory of data envelopment analysis (DEA)—an LP-based performance measurement method proposed by Charnes, Cooper, and Rhodes (1978). See, for example, Charnes, Cooper, and Thrall (1991), Gonzalez-Lima, Tapia, and Thrall (1996), Thompson, Langemeier, and Thrall (1994), Chen, Morita, and Zhu (2003), Sueyoshi and Sekitani (2007a), Sueyoshi and Sekitani (2007b), Sueyoshi and Sekitani (2009), Sueyoshi, Goto (2010), Sueyoshi, Goto (2012a), Sueyoshi and Goto (2012b) and Sueyoshi, Goto (2013).

As demonstrated by Greenberg (1986) and Greenberg (1994), a pair of primal–dual optimal solutions is unique if and only if it is a strictly complementary basic solution. This indicates that a strictly complementary solution can be identified easily in the case of non-degeneracy. The difficulty arises, however, when multiple primal–dual optimal solutions occur due to the degeneracy. To deal with this difficulty, significant research efforts have been devoted. Balinski and Tucker (1969) introduced a so-called Balinski–Tucker simplex tableau to construct a strictly complementary solution. Conversely, Zhang (1994) and Tijssen and Sierksma (1998) showed that a Balinski–Tucker tableau can be generated from a given strictly complementary solution in strongly polynomial time. In some sense, a Balinski–Tucker simplex tableau and a strictly complementary solution are equivalent.

---

[1] This condition is known by several names in the literature as extended complementary slackness (Charnes and Cooper, 1961), full complementary slackness (Balinski, 1968), complete complementary slackness (Balinski and Tucker, 1969), strong complementary slackness (Spivey and Thrall, 1970; Sierksma, 2002), and strict complementary slackness (Tijssen and Sierksma, 1998; Tijssen, 2000).

Beside the interior point methods, which generate a sequence of solutions converging to a strictly complementary solution, the identification issue can be dealt with by resorting to other algorithms adopted for solving LP problems, e.g., the most popular simplex algorithm. In this regard, Pelessoni (1998) proposed an LP problem whose (basic) optimal solution determines a strictly complementary solution. Sierksma and Tijssen (2003) also developed a similar LP problem that was employed later by Sueyoshi and Sekitani (2007a) and Sueyoshi and Sekitani (2007b) for finding all possible references of an observed decision making unit in DEA. However, as Krivonozhko, Førsund, and Lychev (2012b) have argued, not only the computational burden of Sueyoshi and Sekitani's (2007b) approach is high, but it also seems that the basic matrices defined in their approach are likely to be ill-conditioned, leading to erroneous and unacceptable results even for medium-size problems.

The contribution of this study is three fold. Inspired by Motzkin's (1936) idea of making transition from a polyhedral set to a polyhedral cone, we first define a characteristic cone for a given set. We represent this cone as the union of its level sets and show that all of its level sets preserve the structural properties of the given set. We then demonstrate that the convexity of the given set is a necessary and sufficient condition for the convexity of its characteristic cone. Under the convexity and closedness assumptions, we also have some interesting findings.

The first one is that the extreme directions of the characteristic cone are *vertical* and *horizontal*, which are corresponding to extreme points and extreme directions of the given set, respectively. The second finding relates to the characterization of the closure of the characteristic cone. In the case of boundedness, we show that the origin is the only limit point of the characteristic cone that is not in this cone. In the case of unboundedness, however, the conical hull of the horizontal directions is a face of the closure of the characteristic cone whose elements are the only limit points which are not contained in the characteristic cone. The third interesting finding is that the extreme directions of the closure of the characteristic cone are the same vertical and horizontal directions. Based on this finding, we substantiate that the closure of the characteristic cone of a closed convex set containing no line can be decomposed as the direct sum of two convex cones; one generated by the vertical directions, and the other by the horizontal ones. As a main consequence, we then derive the decomposition of a convex set as the direct sum of the convex hull of its extreme points and the conical hull of its extreme directions. For a polyhedron, specially, the decomposition reduces to Motzkin's (1936) decomposition.

As the second contribution, we devise a convex programming (CP) based approach to identify a *maximal element* of a convex set, that is, an element in which the number of positive components is maximum. First, we propose a CP problem to identify a maximal element of a convex cone. Using the proposed CP problem, we then formulate another CP problem to find a maximal element of a convex set. Our formulation is based on the finding that a maximal element of a convex set is corresponding to a maximal element of its characteristic convex cone.

As our last contribution, we propose an LP problem to identify a relative interior point of a polyhedron. It is well established that the relative interior of a polyhedron is exactly the set of its maximal elements. Based on this fact and by using the CP problem developed for finding a maximal element of a convex set, we propose our LP problem. This formulation was motivated by the two facts that the information of a relative interior point of the optimal face of an LP problem is useful for the parametric analysis (Adler and Monteiro, 1992; Mehrotra, 1993), and that each pair of solutions in the relative interiors of the primal and dual optimal faces of an LP problem yields a strictly complementary solution (Nering, 1993). It is worth noting that since our proposed LP problem contains several upper-bounded variables, the computational efficiency of our method can be enhanced by using the simplex algorithm[2] adopted for solving the LP problems with upper-bounded variables, which is much more efficient than the ordinary simplex algorithm (Winston, 2003).

Finally, as an interesting application of our proposed LP problem, we adopt two procedures for identifying a strictly complementary solution. Requiring that the optimal objective value be known, the first procedure determines a pair of relative interior points of the primal and optimal faces by solving two distinct upper-bounded LP problems. In contrast, our second procedure is free from the requirement of the first approach and determines a strictly complementary solution directly by solving a single upper-bounded LP problem.

For solving LP problems, a variety of algorithms have been designed to date, among which the simplex algorithm is very efficient in practice—even though it has exponential worst-case complexity. In Fig.1, we have displayed how the existing algorithms can be employed for yielding a strictly complementary solution for an LP problem. It can be seen that, as well as the methods of Pelessoni (1998) and Sierksma and Tijssen (2003), our proposed method provides the possibility to apply each of the existing algorithms. This is invaluable because a wider availability of alternative methods can facilitate the employment of optimization software.

---

[2] The simplex algorithm for bounded variables was published by Dantzig (1955) and was independently developed by Charnes and Lemke (1954).

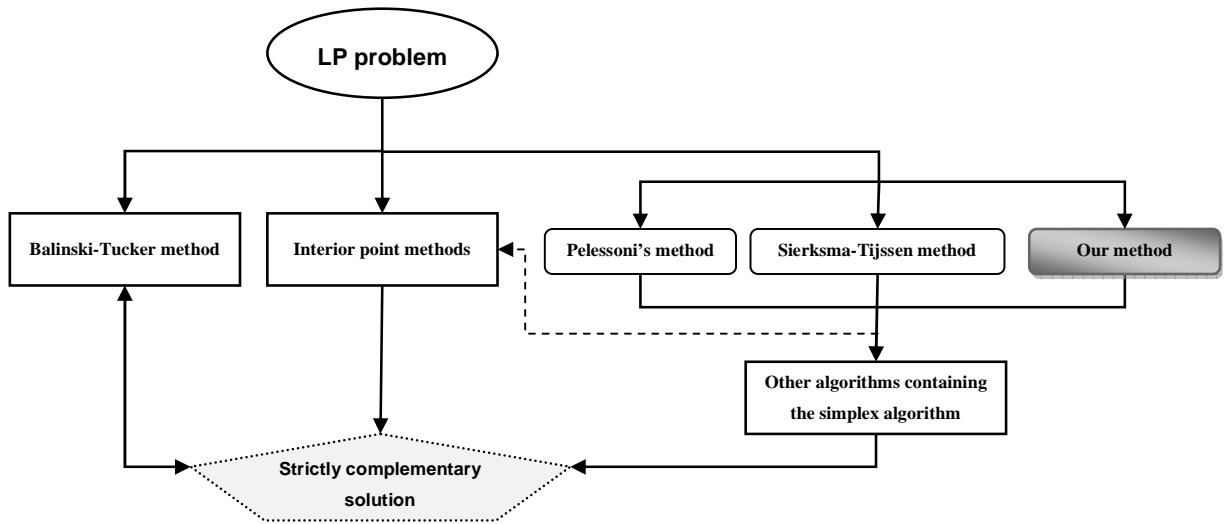

**Fig. 1** Our contribution to the identification of a strictly complementary solution

We note that the problems formulated by Pelessoni (1998) and Sierksma and Tijssen (2003) are, indeed, *max–min* problems. These problems maximize the minimum of the sums of the complementary variables and can be linearized equivalently by a suitable variable change. In comparison with the two above-mentioned methods, our proposed method can be beneficial as its strategy is different depending on whether the optimal value is given. If the optimal value is known, the problem reduces to the solving of two LP problems with small numbers of constrains. Otherwise, an LP problem is solved that has less number of constraints than the linearized forms of the max–min problems developed by Pelessoni (1998) and Sierksma and Tijssen (2003).

Finally, as regards the debate on the preferred choice between our method and the other existing ones for finding a strictly complementary solution, we take no position. We emphasize that the purpose of our *mainly theoretical* study is not to propose a method in the sense of introducing a method in a competition for the fastest running time. Rather, our purpose is to propose an alternative procedure.

The remainder of this paper unfolds as follows. Section 2 discusses the representation of a convex set through its characteristic cone and investigates the properties of this cone. Section 3 proposes a CP problem for finding a maximal element of a convex set. Section 4 presents an LP problem for the identification of a relative interior point of a polyhedron. This problem is then used in Section 5 to develop two methods for finding a strictly complementary solution. Section 6 then illustrates these methods with a numerical example. Section 7 presents the summary of our work with some concluding remarks.

## 2 Characteristic cone of a convex set

As far as notations are concerned, we will use the following notations throughout the paper. We symbolize the sets by capital letters and their members by lower-case letters. We denote vectors and matrices in bold letters, vectors in lower case and matrices in upper case. All vectors are column vectors. We denote by a subscript $T$ the transpose of vectors and matrices. We also use $\mathbf{0}_n$ and $\mathbf{1}_n$ to show $n$-dimensional vectors with the values of 0 and 1 in every entry, respectively. Furthermore, we denote by $\mathbf{I}_n$ the $n \times n$ identity matrix. For $\mathbf{x}, \mathbf{y} \in \mathbb{R}^n$, we write $\mathbf{x} \leq \mathbf{y}$ if $x_j \leq y_j$ for all $j = 1, \ldots, n$.

We begin our discussion with some definitions in convex optimization (Rockafellar, 1970).

**Definition 2.1** A convex subset $F$ of a convex set $X$ in $\mathbb{R}^n$ is called a *face* if and only if, for any $\mathbf{x}_1, \mathbf{x}_2 \in X$ and $\lambda \in (0,1)$, $\lambda \mathbf{x}_1 + (1-\lambda) \mathbf{x}_2 \in F$ implies that $\mathbf{x}_1, \mathbf{x}_2 \in F$. In addition, $F$ is called a *k-face* if its dimension is $k$.

Note that the dimension of a face $F$ equals to that of its corresponding *affine hull*, denoted by $\mathrm{aff}(F)$. Moreover, 0-Faces of $X$ are called its *extreme points*.

**Definition 2.2** A non-zero vector $\mathbf{d}$ in $\mathbb{R}^n$ is called a *recession direction* of a convex set $X$ in $\mathbb{R}^n$ if and only if for each $\mathbf{x} \in X$, $\mathbf{x} + \theta \mathbf{d} \in X$ for all $\theta \geq 0$. The direction $\mathbf{d}$ is *extreme* if and only if $\mathbf{d} = \lambda_1 \mathbf{d}_1 + \lambda_2 \mathbf{d}_2$ with $\lambda_1, \lambda_2 > 0$ implies that $\mathbf{d}_1 = \alpha \mathbf{d}_2$ for some $\alpha > 0$.

77

**Definition 2.3** The *relative interior* of $X$, denoted by $ri(X)$, is defined as the interior of $X$ when this set is considered as a subset of $\text{aff}(X)$. That is, $\mathbf{x}_0 \in ri(X)$ if and only if there exists $\varepsilon > 0$ such that $N_\varepsilon(\mathbf{x}_0) \cap \text{aff}(X) \subseteq X$ [3].

**Definition 2.4** A set $C$ in $\mathbb{R}^n$ is called a *cone* if and only if $\mathbf{x} \in C$ implies that $\alpha \mathbf{x} \in C$ for all $\alpha \geq 0$, whereas a *blunt cone* is a cone without origin, that is, $\mathbf{0}_n \notin C$ and $\mathbf{x} \in C$ implies that $\alpha \mathbf{x} \in C$ for all $\alpha > 0$. Moreover, a *convex cone* is a cone which is convex.

Motzkin (1936) made a transition from a polyhedral set in $\mathbb{R}^n$ to a polyhedral cone in $\mathbb{R}^{n+1}$. Inspired by his innovative idea, we define in general a cone in $\mathbb{R}^{n+1}$ corresponding to a given set $X$ in $\mathbb{R}^n$ as

$$C_X := \left\{ x_{n+1} \begin{pmatrix} \mathbf{x} \\ 1 \end{pmatrix} \mid \mathbf{x} \in X,\ x_{n+1} > 0 \right\}. \tag{1}$$

Obviously, there exists a one-to-one correspondence between $X$ and $C_X$ as

$$\mathbf{x} \in X \Longleftrightarrow \begin{pmatrix} \mathbf{x} \\ 1 \end{pmatrix} \in C_X. \tag{2}$$

Now, we express the relationships between some properties of $X$ and $C_X$ in terms of seven theorems. The following lemma (without proof) states that $C_X$ can be expressed as the union of its level sets.

**Lemma 2.1** For any $\alpha > 0$, let $X^\alpha := \left\{ \alpha \begin{pmatrix} \mathbf{x} \\ 1 \end{pmatrix} \mid \mathbf{x} \in X \right\}$. Then,

(i) $C_X = \bigcup_{\alpha > 0} X^\alpha$.

(ii) $X^\alpha$ is convex for all $\alpha > 0$ if and only if $X$ is convex.

**Theorem 2.1** $C_X$ is a blunt cone. In addition, it is convex if and only if $X$ is convex.

---

[3] Here, $N_\varepsilon(\mathbf{x}_0)$ is the $\varepsilon$-*neighborhood* of $\mathbf{x}_0$ that is defined as $N_\varepsilon(\mathbf{x}_0) := \left\{ \mathbf{x} \mid \|\mathbf{x} - \mathbf{x}_0\| < \varepsilon \right\}$.

See Appendix A for the proof.

**Remark 2.1** For any $\alpha > 0$, let $F^\alpha := \left\{ \alpha \begin{pmatrix} \mathbf{x} \\ 1 \end{pmatrix} \middle| \mathbf{x} \in F \right\}$ where $F$ is a subset of $X$. Then, it can be easily verified that $F^\alpha = X^\alpha \cap C_F$ for all $\alpha > 0$. Moreover, Lemma 2.1 and Theorem 2.1 remain true if $X$ is replaced by $F$.

Now, we propose the following theorem to demonstrate the relationship between $X$ and its $\alpha$-level sets.

**Theorem 2.2** Let $X$ be a convex set in $\mathbb{R}^n$ and let $F$ be a convex subset of $X$. Then,

(i) $\begin{pmatrix} \hat{\mathbf{d}} \\ 0 \end{pmatrix}$ is a recession (an extreme) direction of $X^\alpha$ for all $\alpha > 0$ if and only if $\hat{\mathbf{d}}$ is a recession (an extreme) direction of $X$.

(ii) $F^\alpha$ is a $k$-face (an extreme point) of $X^\alpha$ for all $\alpha > 0$ if and only if $F$ is a $k$-face (extreme point) of $X$.

See Appendix A for the proof.

The following theorem demonstrates the relationship between the faces of $X$ and those of $C_X$.

**Theorem 2.3** Let $X$ be a convex set in $\mathbb{R}^n$ and let $F$ be a convex subset of $X$. Then, $C_F = \bigcup_{\alpha \geq 0} F^\alpha$ is a $(k+1)$-face (extreme direction) of $C_X$ if and only if $F$ is a $k$-face (extreme point) of $X$.

See Appendix A for the proof.

**Theorem 2.4** Let $X$ be a closed convex set in $\mathbb{R}^n$. Then, $\begin{pmatrix} \mathbf{d} \\ d_{n+1} \end{pmatrix}$ is a recession (an extreme) direction of $C_X$ if and only if one of the followings holds:

(i) $d_{n+1} = 0$ and $\mathbf{d}$ is a recession (an extreme) direction of $X$.

(ii) $d_{n+1} > 0$ and $\dfrac{1}{d_{n+1}}\mathbf{d}$ is a (an extreme) point of $X$.

See Appendix A for the proof.

For a closed convex set $X$, Theorem 2.4 states that the extreme directions of $C_X$ are in the two types (i) and (ii), which we call them as *vertical* and *horizontal*, respectively. As the representative sets of vertical and horizontal extreme directions, we also define $D_v$ and $D_h$ as

$$D_v = \left\{ \begin{pmatrix} \mathbf{d} \\ 1 \end{pmatrix} \Big| \; \mathbf{d} \text{ is an extreme point of } X \right\},$$
$$D_h = \left\{ \begin{pmatrix} \mathbf{d} \\ 0 \end{pmatrix} \Big| \; \mathbf{d} \text{ is an extreme direction of } X \text{ with } \|\mathbf{d}\| = 1 \right\}. \tag{3}$$

Now, consider $cl(C_X)$ (the closure of $C_X$) which is the smallest closed convex cone containing $C_X$ (Rockafellar, 1970, page 50). In the cases of unboundedness and boundedness, we characterize the structure of $cl(C_X)$ by showing that the conical hull of $D_h$ ($coni(D_h)$) and the origin ($\{\mathbf{0}_{n+1}\}$) are the respective faces of $cl(C_X)$ whose elements are the only limit points which are not contained in $C_X$.

**Theorem 2.5** Let $X$ be a closed convex set in $\mathbb{R}^n$. Then,
  (i)  $X$ is unbounded if and only if $cl(C_X) = C_X \cup coni(D_h)$.
  (ii) $X$ is bounded if and only if $cl(C_X) = C_X \cup \{\mathbf{0}_{n+1}\}$.

See Appendix A for the proof.

The following theorem states that the extreme directions of $cl(C_X)$ are the same vertical and horizontal directions of $C_X$.

**Theorem 2.6** Let $X$ be a closed convex set in $\mathbb{R}^n$. Then, $D_v \cup D_h$ is the representative set of extreme directions of $cl(C_X)$.

See Appendix A for the proof.

**Corollary 2.1** $cl(C_X)$ is a polyhedral cone if and only if $X$ is a polyhedron.

Now, we present the following theorem to establish that $cl(C_X)$ can be decomposed as the direct sum of two convex cones.

**Theorem 2.7** Let $X$ be a closed convex set in $\mathbb{R}^n$ containing no line. Then,

$$cl(C_X) = cone(conv(D_v)) + coni(D_h). \qquad (4)$$

See Appendix A for the proof.

By setting $z_{n+1} = 1$ in (A.5), we obtain the following important corollary, which states that a convex set can be decomposed as the direct sum of the convex hull of its extreme points and the conical hull of its extreme directions.

**Corollary 2.2** (*Representation theorem*) Let $X$ is a closed convex set in $\mathbb{R}^n$ containing no line and let $E_X$ and $D_X$ be the sets of its extreme points and extreme directions, respectively. Then,

$$X = conv(E_X) + coni(D_X). \qquad (5)$$

Specially, if $X$ is a polyhedron, then (5) yields Motzkin's (1936) decomposition of a polyhedron.

Fig. 2 illustrates the results in the case that $X$ is a polyhedron.

## 3 Finding a maximal element of a convex set

Let $X$ be a convex set in $\mathbb{R}_+^n$. We are concerned with finding a maximal element of $X$ when it is not singleton. To this aim, we present the following definitions.

**Definition 3.1** The *support* of a vector $\mathbf{x} \in X$, $\sigma(\mathbf{x})$, is defined as the index set of its positive components, i.e., $\sigma(\mathbf{x}) = \{j \mid x_j > 0\}$.

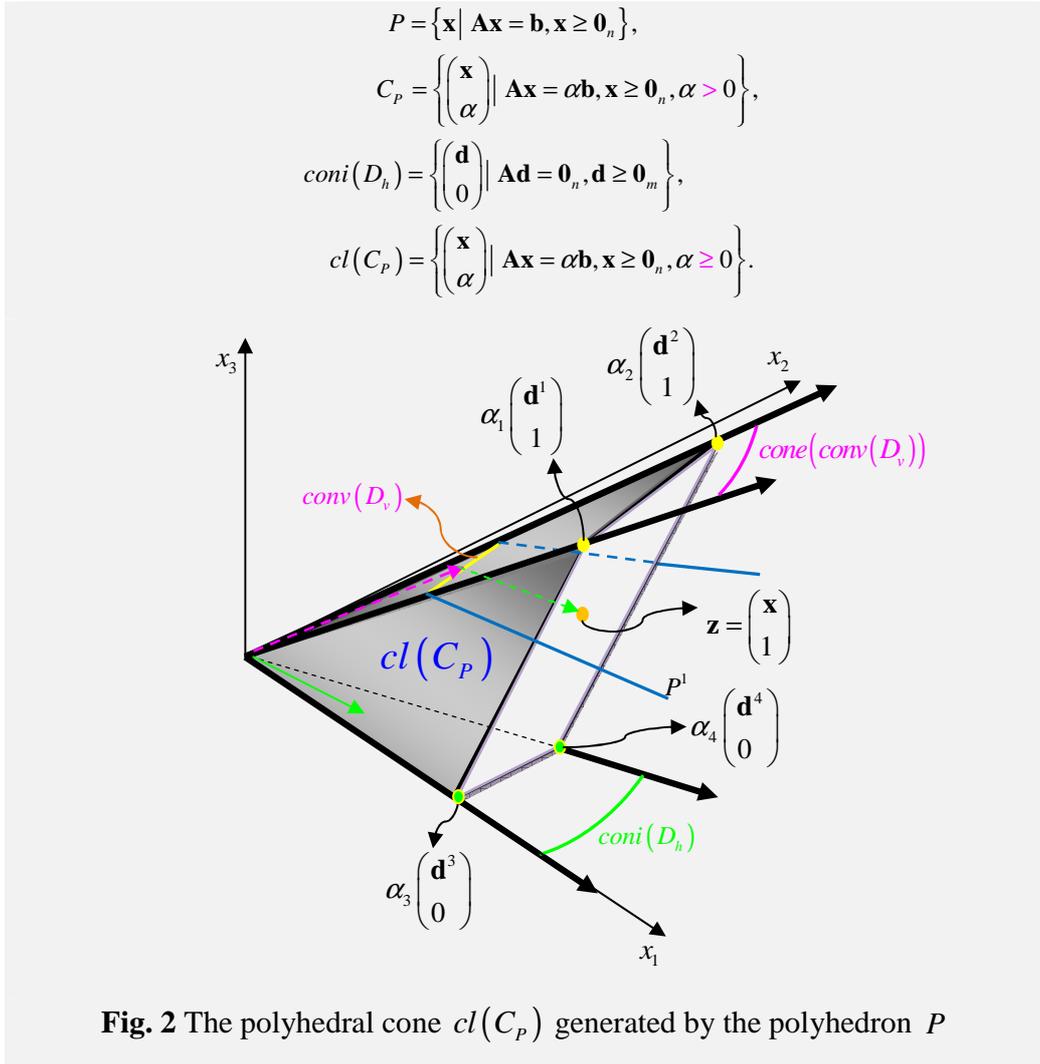

**Fig. 2** The polyhedral cone $cl(C_P)$ generated by the polyhedron $P$

**Definition 3.2** The *zero-norm* of a vector $\mathbf{x} \in X$, $\|\mathbf{x}\|_0$, is defined as the number of its positive components or, equivalently, as the cardinality of its support, i.e., $\|\mathbf{x}\|_0 = Card(\sigma(\mathbf{x}))$.

**Definition 3.3** The *maximal support* of $X$, $\Sigma(X)$, is defined as the support of a *maximal element* of $X$, that is, an element in which the number of positive components is maximum. Formally, $\Sigma(X) = \sigma(\bar{\mathbf{x}})$ where $\bar{\mathbf{x}} \in \arg\max_{\mathbf{x} \in X} \|\mathbf{x}\|_0$.

Now, we state a result concerning the maximal elements of $X$.

**Lemma 3.1** The following statements hold:

(i) $\sigma(\lambda \mathbf{x}_1 + (1-\lambda)\mathbf{x}_2) = \sigma(\mathbf{x}_1) \cup \sigma(\mathbf{x}_2)$ for any $\mathbf{x}_1, \mathbf{x}_2 \in X$ and any $\lambda \in (0,1)$.

(ii) $\Sigma(X) = \bigcup_{\mathbf{x} \in X} \sigma(\mathbf{x})$.

(iii) $\arg\max_{\mathbf{x} \in X} \|\mathbf{x}\|_0 = \{\mathbf{x} \in X \mid x_j > 0, j \in \{1,...,n\} \setminus J^=\}$, where $J^=$ denotes the index set of components which take zero values throughout $X$, i.e., $J^= = \{j \mid x_j = 0 \text{ for all } \mathbf{x} \in X\}$.

(iv) $\bar{\mathbf{x}} \in \arg\max_{\mathbf{x} \in X} \|\mathbf{x}\|_0$ if and only if $\begin{pmatrix}\bar{\mathbf{x}}\\1\end{pmatrix} \in \arg\max_{\binom{\mathbf{x}}{1} \in C_X} \left\|\begin{pmatrix}\mathbf{x}\\1\end{pmatrix}\right\|_0$.

See Appendix A for the proof.

Part (ii) of the above lemma states that $\Sigma(X)$ is unique. Hence, the maximal elements of $X$ are all equivalent in the sense that the indices of positive components are invariant among all of these elements. Moreover, $\Sigma(X)$ contains the support of all $\mathbf{x} \in X$. Based on this fact, we develop a simple algorithm for finding a maximal element of $X$ as follows:

**Step 1** (*Initialization*) Start with $J_1 = \{1,...,n\}$ and $I_1 = \emptyset$. Set $k = 1$.

**Step 2** Solve the following convex programming (CP) problem:

$$[\text{CP}_k] \quad \max \sum_{j \in J_k} x_j$$
$$\text{s.t.} \quad \mathbf{x} \in X.$$

**Step 3** Let $\mathbf{x}^k$ be an optimal solution for $\text{CP}_k$. Then, set $J_{k+1} := J_k - \sigma(\mathbf{x}^k)$ and $I_{k+1} := I_k \cup \sigma(\mathbf{x}^k)$.

**Step 4** If $J_{k+1} \subset J_k$, set $k \leftarrow k+1$ and return to Step 2; otherwise, the algorithm is terminated.

If the algorithm terminates at iteration $K$, then it returns the maximal support and a maximal element of $X$ as $\Sigma(X) = I_K$ and $\bar{\mathbf{x}} = \frac{1}{K}\sum_{k=1}^{K}\mathbf{x}^k$, respectively.

While the implementation of the above algorithm is straightforward, it may require solving many problems when $n$ is considerably large. This motivates us to propose a new approach that needs the execution of a single problem. In order to develop our approach, we begin by the following fundamental theorem which determines a maximal element of a convex cone.

**Theorem 3.1** Let $C$ be a convex cone in $\mathbb{R}_+^{n+1}$ and let $(\mathbf{z}^{1*}, \mathbf{z}^{2*})$ be an optimal solution to the following CP problem

$$\begin{aligned} \max \quad & \sum_{j=1}^{n+1} z_j^2 \\ \text{s.t.} \quad & \mathbf{z}^1 + \mathbf{z}^2 \in C, \\ & \mathbf{z}^1 \geq \mathbf{0}_{n+1}, \ \mathbf{z}^2 \leq \mathbf{1}_{n+1}. \end{aligned} \quad (6)$$

Then,

(i) $\mathbf{z}^{2*} \geq \mathbf{0}_{n+1}$.

(ii) $z_j^{2*} = 1$ for every $j$ that $z_j^{2*} > 0$.

(iii) $\overline{\mathbf{z}} := \mathbf{z}^{1*} + \mathbf{z}^{2*} \in \arg\max_{\mathbf{z} \in C} \|\mathbf{z}\|_0$.

See Appendix A for the proof.

As mentioned earlier, our aim is to find a maximal element of the convex set $X$. Based on part (iv) of Lemma 3.1, a maximal element of $X$ can be found by identifying a maximal element of $C_X$. This fact enables us to effectively apply Theorem 3.1 for our aim as established in the following theorem.

**Theorem 3.2** Let $X$ be a convex set in $\mathbb{R}_+^n$ and let $\left(\begin{pmatrix} \mathbf{x}^{1*} \\ w^{1*} \end{pmatrix}, \begin{pmatrix} \mathbf{x}^{2*} \\ w^{2*} \end{pmatrix}\right)$ be an optimal solution to the following CP problem:

$$\begin{aligned} \max \quad & \sum_{j=1}^{n} x_j^2 + w^2 \\ \text{s.t.} \quad & \begin{pmatrix} \mathbf{x}^1 \\ w^1 \end{pmatrix} + \begin{pmatrix} \mathbf{x}^2 \\ w^2 \end{pmatrix} \in C_X, \\ & \begin{pmatrix} \mathbf{x}^1 \\ w^1 \end{pmatrix} \geq \mathbf{0}_{n+1}, \ \begin{pmatrix} \mathbf{x}^2 \\ w^2 \end{pmatrix} \leq \mathbf{1}_{n+1}. \end{aligned} \quad (7)$$

If $w^{2*} = 0$, then $X = \emptyset$. Otherwise, $w^{2*} = 1$ and $\bar{\mathbf{x}} := \dfrac{1}{w^{1*}+1}\left(\mathbf{x}^{1*} + \mathbf{x}^{2*}\right)$ lies in $\arg\max\limits_{\mathbf{x} \in X} \|\mathbf{x}\|_0$.

See Appendix A for the proof.

**Remark 3.1** Note that $cl(C_X)$ can be used also instead of $C_X$ in model (7). This is because maximality occurs always at $C_X$.

## 4 Finding a relative interior point of a polyhedron

Let $P$ be a polyhedral set defined by

$$P = \{\mathbf{x} \mid \mathbf{A}\mathbf{x} = \mathbf{b}, \mathbf{x} \geq \mathbf{0}_n\}, \tag{8}$$

where $\mathbf{A}$ is a matrix of order $m \times n$ and $\mathbf{b} \in \mathbb{R}^m$.

We aim to find a relative interior point of $P$. To this aim, we first establish an important relationship between relative interior points and maximal elements of $P$.

**Theorem 4.1** $ri(P) = \arg\max\limits_{\mathbf{x} \in P} \|\mathbf{x}\|_0$.

See Appendix A for the proof.

The above theorem states that the relative interior of $P$ consists of its maximal elements. This follows that the problem of finding a relative interior point of $P$ is equivalent to finding a maximal element of this set. Hence, as a useful result of Theorems 3.2 and 4.1, we have the following corollary.

**Corollary 4.1** Let $\left(\begin{pmatrix}\mathbf{x}^{1*}\\w^{1*}\end{pmatrix}, \begin{pmatrix}\mathbf{x}^{2*}\\w^{2*}\end{pmatrix}\right)$ be an optimal solution to the following LP problem:

$$\begin{aligned}
\max \quad & \sum_{j=1}^{n} x_j^2 + w^2 \\
\text{s.t.} \quad & [\mathbf{A}, -\mathbf{b}]\begin{bmatrix}\mathbf{x}^1 + \mathbf{x}^2 \\ w^1 + w^2\end{bmatrix} = \mathbf{0}_m, \\
& \begin{pmatrix}\mathbf{x}^1 \\ w^1\end{pmatrix} \geq \mathbf{0}_{n+1}, \quad \mathbf{0}_{n+1} \leq \begin{pmatrix}\mathbf{x}^2 \\ w^2\end{pmatrix} \leq \mathbf{1}_{n+1}.
\end{aligned} \tag{9}$$

If $w^{2*} = 0$, then $P = \emptyset$. Otherwise, $w^{2*} = 1$ and $\overline{\mathbf{x}} := \dfrac{1}{w^{1*}+1}\left(\mathbf{x}^{1*} + \mathbf{x}^{2*}\right) \in ri(P)$.

**Remark 4.1** The constraints $\mathbf{0}_{n+1} \leq \begin{pmatrix} \mathbf{x}^2 \\ w^2 \end{pmatrix}$ are added to (9) in order to guarantee that $\begin{bmatrix} \mathbf{x}^1 + \mathbf{x}^2 \\ w^1 + w^2 \end{bmatrix} \in cl(C_P)$.

It is worth noting that since model (9) contains several upper-bounded variables, the computational efficiency can be enhanced by using the simplex algorithm adopted for solving the LP problems with upper-bounded variables, which is much more efficient than the ordinary simplex algorithm (Winston, 2003). Indeed, if model (9) is considered as an LP problem with upper-bounded variables, then the number of its constraints reduces from $m+n+1$ to $m$.

## 5 Finding a strictly complementary solution

In this section, we employ model (9) for identifying a strictly complementary solution. First, consider the primal LP problem in canonical form

$$(P): \quad \max_{\mathbf{x} \in X} \mathbf{c}^T \mathbf{x} \quad \text{subject to} \quad X = \{\mathbf{x} \mid \mathbf{A}\mathbf{x} \leq \mathbf{b},\ \mathbf{x} \geq \mathbf{0}\}, \qquad (10)$$

and its dual

$$(D): \quad \min_{\mathbf{y} \in Y} \mathbf{b}^T \mathbf{y} \quad \text{subject to} \quad Y = \{\mathbf{y} \mid \mathbf{A}^T \mathbf{y} \geq \mathbf{c},\ \mathbf{y} \geq \mathbf{0}\}, \qquad (11)$$

where $\mathbf{c} \in \mathbb{R}^n$, $\mathbf{b} \in \mathbb{R}^m$ are fixed data and $\mathbf{x} \in \mathbb{R}^n$ and $\mathbf{y} \in \mathbb{R}^m$ are unknown vectors.

By including variables $\mathbf{v} \in \mathbb{R}^m_+$ and $\mathbf{u} \in \mathbb{R}^n_+$ into the inequality constraints of problems (P) and (D), we obtain

$$(P): \quad \max_{\binom{\mathbf{x}}{\mathbf{u}} \in S_P} \mathbf{c}^T \mathbf{x} \quad \text{subject to} \quad S_P = \left\{ \begin{pmatrix} \mathbf{x} \\ \mathbf{u} \end{pmatrix} \middle| [\mathbf{A}, \mathbf{I}_m] \begin{bmatrix} \mathbf{x} \\ \mathbf{u} \end{bmatrix} = \mathbf{b},\ \begin{pmatrix} \mathbf{x} \\ \mathbf{u} \end{pmatrix} \geq \mathbf{0}_{n+m} \right\},$$

$$(D): \quad \min_{\binom{\mathbf{y}}{\mathbf{v}} \in S_D} \mathbf{b}^T \mathbf{y} \quad \text{subject to} \quad S_D = \left\{ \begin{pmatrix} \mathbf{y} \\ \mathbf{v} \end{pmatrix} \middle| [\mathbf{A}^T, -\mathbf{I}_n] \begin{bmatrix} \mathbf{y} \\ \mathbf{v} \end{bmatrix} = \mathbf{c},\ \begin{pmatrix} \mathbf{y} \\ \mathbf{v} \end{pmatrix} \geq \mathbf{0}_{m+n} \right\}, \qquad (12)$$

where $\mathbf{I}_m$ and $\mathbf{I}_n$ are identity matrices of appropriate sizes.

Any vectors $\begin{pmatrix} \mathbf{x} \\ \mathbf{u} \end{pmatrix} \in S_P$ and $\begin{pmatrix} \mathbf{y} \\ \mathbf{v} \end{pmatrix} \in S_D$ are called primal and dual feasible solutions of (P) and (D), respectively. Feasible solutions $\begin{pmatrix} \mathbf{x} \\ \mathbf{u} \end{pmatrix} \in S_P$ and $\begin{pmatrix} \mathbf{y} \\ \mathbf{v} \end{pmatrix} \in S_D$ to (P) and (D) are optimal if and only if the CSC holds, i.e., $\mathbf{x}^T \mathbf{v} = 0$ and $\mathbf{y}^T \mathbf{u} = 0$. Then, the pairs $(x_j, v_j)$, $j=1,...,n$, and $(y_i, u_i)$, $i=1,...,m$, are called *complementary variables*.

Assuming that problems (P) and (D) have finite optimal objective value, the SCSC implies the existence of optimal solutions $\begin{pmatrix} \mathbf{x}^* \\ \mathbf{u}^* \end{pmatrix} \in S_P$ and $\begin{pmatrix} \mathbf{y}^* \\ \mathbf{v}^* \end{pmatrix} \in S_D$ to (P) and (D) for which $\mathbf{x}^* + \mathbf{v}^* > \mathbf{0}_n$ and $\mathbf{y}^* + \mathbf{u}^* > \mathbf{0}_m$. This result is known as *Goldman-Tucker theorem* (Goldman and Tucker, 1956) and the pair $\left( \begin{pmatrix} \mathbf{x}^* \\ \mathbf{u}^* \end{pmatrix}, \begin{pmatrix} \mathbf{y}^* \\ \mathbf{v}^* \end{pmatrix} \right)$ is called a *strictly complementary solution*. This theorem states that, for each pair of complementary variables, exactly one of the variables takes positive value and the other is zero. Moreover, the following unique partitions, referred to as *optimal partitions*, can be determined for the index sets $\{1,...,n\}$ and $\{1,...,m\}$:

$$\sigma(\mathbf{x}^*) \cup \sigma(\mathbf{v}^*) = \{1,...,n\}, \quad \sigma(\mathbf{x}^*) \cap \sigma(\mathbf{v}^*) = \emptyset,$$
$$\sigma(\mathbf{u}^*) \cup \sigma(\mathbf{y}^*) = \{1,...,m\}, \quad \sigma(\mathbf{u}^*) \cap \sigma(\mathbf{y}^*) = \emptyset. \tag{13}$$

Let $F_P^*$ and $F_D^*$ denote the optimal faces of problems (P) and (D), respectively. Then, the relationship between the two optimal faces $F_P^*$ and $F_D^*$ and the SCSC is demonstrated by the following theorem (Nering and Tucker, 1993; Zhang, 1994).

**Theorem 5.1** A solution pair $\left( \begin{pmatrix} \mathbf{x}^* \\ \mathbf{u}^* \end{pmatrix}, \begin{pmatrix} \mathbf{y}^* \\ \mathbf{v}^* \end{pmatrix} \right)$ is a strictly complementary solution if and only if $\begin{pmatrix} \mathbf{x}^* \\ \mathbf{u}^* \end{pmatrix} \in ri(F_P^*)$ and $\begin{pmatrix} \mathbf{y}^* \\ \mathbf{v}^* \end{pmatrix} \in ri(F_D^*)$.

### *5.1 First approach*

This approach is based on the primary assumption that the optimal objectively value of the primal and dual LP problems (10) and (11) is known. Let $z^*$ be the optimal objective value for these problems. Then, the hyperplanes $H_P^* = \left\{ \begin{pmatrix} \mathbf{x} \\ \mathbf{u} \end{pmatrix} \mid \mathbf{c}^T \mathbf{x} = z^* \right\}$ and $H_D^* = \left\{ \begin{pmatrix} \mathbf{y} \\ \mathbf{v} \end{pmatrix} \mid \mathbf{b}^T \mathbf{y} = z^* \right\}$ support $S_P$ and $S_D$, respectively. Consequently, $F_P^*$ and $F_D^*$ can be formulated as

$$
\begin{aligned}
F_P^* &= H_P^* \cap S_P = \left\{ \begin{pmatrix} \mathbf{x} \\ \mathbf{u} \end{pmatrix} \Bigg| \begin{bmatrix} \mathbf{A} & \mathbf{I}_m \\ \mathbf{c}^T & \mathbf{0}_m^T \end{bmatrix} \begin{bmatrix} \mathbf{x} \\ \mathbf{u} \end{bmatrix} = \begin{bmatrix} \mathbf{b} \\ z^* \end{bmatrix}, \begin{pmatrix} \mathbf{x} \\ \mathbf{u} \end{pmatrix} \geq \mathbf{0}_{n+m} \right\}, \\
F_D^* &= H_D^* \cap S_D = \left\{ \begin{pmatrix} \mathbf{y} \\ \mathbf{v} \end{pmatrix} \Bigg| \begin{bmatrix} \mathbf{A}^T & -\mathbf{I}_n \\ \mathbf{b}^T & \mathbf{0}_n^T \end{bmatrix} \begin{bmatrix} \mathbf{y} \\ \mathbf{v} \end{bmatrix} = \begin{bmatrix} \mathbf{c} \\ z^* \end{bmatrix}, \begin{pmatrix} \mathbf{y} \\ \mathbf{v} \end{pmatrix} \geq \mathbf{0}_{m+n} \right\}.
\end{aligned}
\qquad (14)
$$

Theorem 5.1 states that the SCSC is valid on the relative interiors of $F_P^*$ and $F_D^*$. Therefore, applying Corollary 4.1, we propose the following two LP problems to find a strictly complementary solution:

$$
\begin{aligned}
\max \quad & \sum_{j=1}^{n} x_j^2 + \sum_{i=1}^{m} u_i^2 + w_P^2 \\
\text{s.t.} \quad & \begin{bmatrix} \mathbf{A} & \mathbf{I}_m & -\mathbf{b} \\ \mathbf{c}^T & \mathbf{0}_m^T & -z^* \end{bmatrix} \begin{bmatrix} \mathbf{x}^1 + \mathbf{x}^2 \\ \mathbf{u}^1 + \mathbf{u}^2 \\ w_P^1 + w_P^2 \end{bmatrix} = \mathbf{0}_{m+1}, \\
& \begin{pmatrix} \mathbf{x}^1 \\ \mathbf{u}^1 \\ w_P^1 \end{pmatrix} \geq \mathbf{0}_{n+m+1}, \ \mathbf{0}_{n+m+1} \leq \begin{pmatrix} \mathbf{x}^2 \\ \mathbf{u}^2 \\ w_P^2 \end{pmatrix} \leq \mathbf{1}_{n+m+1}.
\end{aligned}
\qquad (15)
$$

and

$$
\begin{aligned}
\max \quad & \sum_{i=1}^{m} y_i^2 + \sum_{j=1}^{n} v_j^2 + w_D^2 \\
\text{s.t.} \quad & \begin{bmatrix} \mathbf{A}^T & -\mathbf{I}_n & -\mathbf{c} \\ \mathbf{b}^T & \mathbf{0}_n^T & -z^* \end{bmatrix} \begin{bmatrix} \mathbf{y}^1 + \mathbf{y}^2 \\ \mathbf{v}^1 + \mathbf{v}^2 \\ w_D^1 + w_D^2 \end{bmatrix} = \mathbf{0}_{n+1}, \\
& \begin{pmatrix} \mathbf{y}^1 \\ \mathbf{v}^1 \\ w_D^1 \end{pmatrix} \geq \mathbf{0}_{m+n+1}, \ \mathbf{0}_{m+n+1} \leq \begin{pmatrix} \mathbf{y}^2 \\ \mathbf{v}^2 \\ w_D^2 \end{pmatrix} \leq \mathbf{1}_{m+n+1}.
\end{aligned}
\qquad (16)
$$

Let $\left( \begin{pmatrix} \mathbf{x}^{1*} \\ \mathbf{u}^{1*} \\ w_P^{1*} \end{pmatrix}, \begin{pmatrix} \mathbf{x}^{2*} \\ \mathbf{u}^{2*} \\ w_P^{1*} \end{pmatrix} \right)$ and $\left( \begin{pmatrix} \mathbf{y}^{1*} \\ \mathbf{v}^{1*} \\ w_D^{2*} \end{pmatrix}, \begin{pmatrix} \mathbf{y}^{2*} \\ \mathbf{v}^{2*} \\ w_D^{2*} \end{pmatrix} \right)$ be optimal solutions to (15) and (16), respectively.

Then, we have

$$\begin{pmatrix} \overline{\mathbf{x}} \\ \overline{\mathbf{u}} \end{pmatrix} := \frac{1}{w_P^{1*}+1} \begin{pmatrix} \mathbf{x}^{1*} + \mathbf{x}^{2*} \\ \mathbf{u}^{1*} + \mathbf{u}^{2*} \end{pmatrix} \in ri\left(F_P^*\right), \tag{17}$$

and

$$\begin{pmatrix} \overline{\mathbf{y}} \\ \overline{\mathbf{v}} \end{pmatrix} := \frac{1}{w_D^{1*}+1} \begin{pmatrix} \mathbf{y}^{1*} + \mathbf{y}^{2*} \\ \mathbf{v}^{1*} + \mathbf{v}^{2*} \end{pmatrix} \in ri\left(F_D^*\right). \tag{18}$$

Therefore, $\left( \begin{pmatrix} \overline{\mathbf{x}} \\ \overline{\mathbf{u}} \end{pmatrix}, \begin{pmatrix} \overline{\mathbf{y}} \\ \overline{\mathbf{v}} \end{pmatrix} \right)$ is a strictly complementary solution.

**Remark 5.1.1** If $n$ is considerably larger that $m$ ($n \gg m$), then solving of problem (16) will be difficult as the number of constraints increases. To enhance the computational efficiency in such a case, we suggest two techniques. The first one is to solve the dual of this problem more quickly and correctly by using the generalized upper-bounding techniques (See Dantzig and Van Slyke (1967).). In the second technique, problem (15) is first solved to obtain $\begin{pmatrix} \overline{\mathbf{x}} \\ \overline{\mathbf{u}} \end{pmatrix} \in ri\left(F_P^*\right)$ as in (17).

Then, based on (13) and (14), the dual of the following problem (after linearization) is solved:

$$\begin{aligned} \max \quad & \min_{i \notin \sigma(\overline{u}_i), j \notin \sigma(\overline{x}_j)} \{y_i, v_j\} \\ \text{s.t.} \quad & \begin{bmatrix} \mathbf{A}^T & -\mathbf{I}_n \\ \mathbf{b}^T & \mathbf{0}_n^T \end{bmatrix} \begin{bmatrix} \mathbf{y} \\ \mathbf{v} \end{bmatrix} = \begin{bmatrix} \mathbf{c} \\ z^* \end{bmatrix}, \\ & \begin{pmatrix} \mathbf{y} \\ \mathbf{v} \end{pmatrix} \geq \mathbf{0}_{m+n}. \end{aligned} \tag{19}$$

### 5.2 Second approach

Consider the set of all pairs of the primal–dual optimal solutions, $F_{PD}^*$, formulated as follows:

$$F_{PD}^* = F_P^* \times F_D^* = \left\{ \begin{pmatrix} \mathbf{x} \\ \mathbf{u} \\ \mathbf{y} \\ \mathbf{v} \end{pmatrix} \middle| \begin{bmatrix} \mathbf{A} & \mathbf{I}_m & \mathbf{0}_{m \times m} & \mathbf{0}_{m \times n} \\ \mathbf{0}_{n \times n} & \mathbf{0}_{n \times m} & \mathbf{A}^T & -\mathbf{I}_n \\ \mathbf{c}^T & \mathbf{0}_m^T & -\mathbf{b}^T & \mathbf{0}_n^T \end{bmatrix} \begin{bmatrix} \mathbf{x} \\ \mathbf{u} \\ \mathbf{y} \\ \mathbf{v} \end{bmatrix} = \begin{bmatrix} \mathbf{b} \\ \mathbf{c} \\ 0 \end{bmatrix}, \begin{bmatrix} \mathbf{x} \\ \mathbf{u} \\ \mathbf{y} \\ \mathbf{v} \end{bmatrix} \geq \mathbf{0}_{2m+2n} \right\}. \quad (20)$$

The set of constraints defining $F_{PD}^*$ are the sets of constraints of problems (P) and (D) together with the constraint $\mathbf{c}^T\mathbf{x} - \mathbf{b}^T\mathbf{y} = 0$, which equalizes the primal and dual objective values to ensure that $\begin{pmatrix} \mathbf{x} \\ \mathbf{u} \end{pmatrix} \in F_P^*$ and $\begin{pmatrix} \mathbf{y} \\ \mathbf{v} \end{pmatrix} \in F_D^*$.

To identify a strictly complementary solution, we find a maximal element of $F_{PD}^*$. To this purpose, we develop the following LP problem by Corollary 4.1:

$$\max \quad \sum_{j=1}^n x_j^2 + \sum_{i=1}^m u_i^2 + \sum_{i=1}^m y_i^2 + \sum_{j=1}^n v_j^2 + w^2$$

$$\text{s.t.} \quad \begin{bmatrix} \mathbf{A} & \mathbf{I}_m & \mathbf{0}_{m \times m} & \mathbf{0}_{m \times n} & -\mathbf{b} \\ \mathbf{0}_{n \times n} & \mathbf{0}_{n \times m} & \mathbf{A}^T & -\mathbf{I}_n & -\mathbf{c} \\ \mathbf{c}^T & \mathbf{0}_m^T & -\mathbf{b}^T & \mathbf{0}_n^T & 0 \end{bmatrix} \begin{bmatrix} \mathbf{x}^1 + \mathbf{x}^2 \\ \mathbf{u}^1 + \mathbf{u}^2 \\ \mathbf{y}^1 + \mathbf{y}^2 \\ \mathbf{v}^1 + \mathbf{v}^2 \\ w^1 + w^2 \end{bmatrix} = \mathbf{0}_{m+n+1},$$

$$\begin{bmatrix} \mathbf{x}^1 \\ \mathbf{u}^1 \\ \mathbf{y}^1 \\ \mathbf{v}^1 \\ w^1 \end{bmatrix} \geq \mathbf{0}_{2m+2n+1}, \quad \mathbf{0}_{2m+2n+1} \leq \begin{bmatrix} \mathbf{x}^2 \\ \mathbf{u}^2 \\ \mathbf{y}^2 \\ \mathbf{v}^2 \\ w^2 \end{bmatrix} \leq \mathbf{1}_{2m+2n+1}.$$

$$(21)$$

Let $\left( \begin{pmatrix} \mathbf{x}^{1*} \\ \mathbf{u}^{1*} \\ \mathbf{y}^{1*} \\ \mathbf{v}^{1*} \\ w^{1*} \end{pmatrix}, \begin{pmatrix} \mathbf{x}^{2*} \\ \mathbf{u}^{2*} \\ \mathbf{y}^{2*} \\ \mathbf{v}^{2*} \\ w^{2*} \end{pmatrix} \right)$ be an optimal solutions to (21). Then, there are two cases: $w^{2*} = 0$ or $w^{2*} = 1$. The first case follows that $F_{PD}^* = \varnothing$; or, equivalently, $F_P^* = \varnothing$ or $F_D^* = \varnothing$. In the second case, we have

$$\begin{pmatrix} \bar{\mathbf{x}} \\ \bar{\mathbf{u}} \end{pmatrix} := \frac{1}{w^{1*}+1} \begin{pmatrix} \mathbf{x}^{1*} + \mathbf{x}^{2*} \\ \mathbf{u}^{1*} + \mathbf{u}^{2*} \end{pmatrix} \in ri(F_P^*), \tag{22}$$

and

$$\begin{pmatrix} \bar{\mathbf{y}} \\ \bar{\mathbf{v}} \end{pmatrix} := \frac{1}{w^{1*}+1} \begin{pmatrix} \mathbf{y}^{1*} + \mathbf{y}^{2*} \\ \mathbf{v}^{1*} + \mathbf{v}^{2*} \end{pmatrix} \in ri(F_D^*), \tag{23}$$

indicating that $\left( \begin{pmatrix} \bar{\mathbf{x}} \\ \bar{\mathbf{u}} \end{pmatrix}, \begin{pmatrix} \bar{\mathbf{y}} \\ \bar{\mathbf{v}} \end{pmatrix} \right)$ is a strictly complementary solution.

## 6 Numerical example

We illustrate our proposed method with the following LP problems, which are taken from Tijssen and Sierksma (1998). To carry out all the computations in this example, we have developed a computer program, given in Appendix B, using the GAMS (General Algebraic Modeling System) optimization software.

Consider the LP problem

$$\begin{aligned}
\max \quad & -4x_1 + 4x_2 - 8x_3 + 4x_4 \\
\text{s.t.} \quad & -1x_1 + 1x_2 - 2x_3 + 1x_4 + u_1 &&= 1, \\
& +4x_1 - 4x_2 + 1x_3 - 2x_4 + u_2 &&= 0, \\
& \phantom{+4x_1 - 4x_2 +\ } -3x_3 + 1x_4 + u_3 &&= 2, \\
& -1x_1 + 1x_2 - 2x_3 + 1x_4 + u_4 &&= 1, \\
& -2x_1 + 5x_2 - 9x_3 + 3x_4 + u_5 &&= 7, \\
& x_j \geq 0,\ j=1,\ldots,4,\ u_i \geq 0,\ i=1,\ldots,5.
\end{aligned} \tag{24}$$

and its corresponding dual problem

$$\begin{aligned}
\min \quad & +1y_1 + 2y_3 + 1y_4 + 7y_5 \\
\text{s.t.} \quad & -1y_1 + 4y_2 \phantom{{}+1y_3} -1y_4 - 2y_5 - v_1 &&= -4, \\
& +1y_1 - 4y_2 \phantom{{}+1y_3} +1y_4 + 5y_5 - v_2 &&= 4, \\
& -2y_1 + 1y_2 - 3y_3 - 2y_4 - 9y_5 - v_3 &&= -8, \\
& +1y_1 - 2y_2 + 1y_3 + 1y_4 + 3y_5 - v_4 &&= 4, \\
& y_i \geq 0,\ i=1,\ldots,5,\ v_j \geq 0,\ j=1,\ldots,4.
\end{aligned} \tag{25}$$

Note that the original constraints of the above LP problems were stated as inequalities that are now converted to equations by adding slack variables $u_i$, $i=1,\ldots,5$ and $v_j$, $j=1,\ldots 4$.

Solving one of problems (24) and (25) obtains the optimal objective value $z^* = 4$. Having this value, we solve models (15) and (16) to find the strictly complementary solution

$$x_1^* = 2.667, \ x_2^* = 1.667, \ x_3^* = 1, \ x_4^* = 4, \ u_1^* = 0, \ u_2^* = 3, \ u_3^* = 1, \ u_4^* = 0, \ u_5^* = 1;$$
$$v_1^* = 0, \quad v_2^* = 0, \quad v_3^* = 0, \ v_4^* = 0, \ y_1^* = 3, \ y_2^* = 0, \ y_3^* = 0, \ y_4^* = 1, \ y_5^* = 0. \quad (26)$$

From (26), $\sigma(\mathbf{x}^*) = \{1,2,3,4\}$ and $\sigma(\mathbf{v}^*) = \varnothing$ constitute the optimal partition of the set $\{1,...,4\}$ ($n = 4$). Further, $\sigma(\mathbf{u}^*) = \{2,3,5\}$ and $\sigma(\mathbf{y}^*) = \{1,4\}$ yield the optimal partition of the set $\{1,...,5\}$ ($m = 5$). Note that the results obtained from model (21) are the same as those presented above.

Since the variables taking zero value in (26) are zero over the corresponding optimal faces, the first and fourth constraints of (24) and all the four constraints of (25) are binding over the primal and dual optimal faces, respectively.

## 7 Conclusion

The present study was primarily aimed at the identification of a maximal element of a convex set and at its application for finding a strictly complementary solution in linear programming. This application is important due to the fact that the identification of a strictly complementary solution is difficult when multiple optimal solutions occur due to the degeneracy.

Based on the proof of Motzkin's (1936) decomposition theorem, we first represented a convex set by its corresponding characteristic cone. Then, we investigated properties of the characteristic cone in details. In this investigation, we demonstrated that the closure of the characteristic cone of a closed convex set containing no line can be decomposed as the direct sum of two convex cones; one generated by the vertical directions, and the other by the horizontal ones. Note that the vertical and horizontal directions are the extreme directions of the characteristic cone, which are corresponding to extreme points and extreme directions of the given set, respectively. As a secondarily objective, we next arrived at an important conclusion that a closed convex set containing no line can be expressed as the direct sum of the convex hull of its extreme points and conical hull of its extreme directions. This result is, indeed, the extension of Motzkin's (1936) decomposition theorem for the convex sets.

In order to find a maximal element of a convex set, we showed that the maximal elements of a convex set and those of its characteristic cone are corresponding. Based on this fact, we equivalently turned the problem under consideration into finding a maximal element of the characteristic cone. To address this problem, we first developed a CP problem to find a maximal element of a convex cone. Using this CP problem, we then formulated another CP problem which determines a maximal element of the characteristic cone and whereby a maximal element of the given convex set.

It was established that the relative interior of a polyhedron is characterized by its maximal elements. In view of this fact and by using the CP problem determining a maximal element of the characteristic cone, we formulated an LP problem to identify a relative interior point of a polyhedron, as our tertiary objective. Then, using the proposed LP problem, we developed two approaches for finding a strictly complementary solution. Knowing the optimal objective value, the optimal faces of the primal and dual LP problems are formulated in the first approach. Then, two distinct LP problems are solved over the primal and dual optimal faces to identify a relative interior point for each of these two faces. The second approach is free from the requirement of the first approach and accomplishes the task by solving of another LP problem. Comparing the two approaches, each has advantages and disadvantages over each other. While the first approach requires knowing the optimal objective value and solving of two LP problems as well, it can be efficient than the second one from the solvability point of view. This is because the numbers of constraints of the two LP problems solved in the first approach are both less than that of the single LP problem of the second approach.

To date, a variety of algorithms have been developed for solving LP problems among which the simplex algorithm is very popular and very efficient in practice—even though it has exponential worst-case complexity. Our proposed two approaches provide the possibility to apply any of the algorithms adopted for solving LP problems including the simplex algorithm. Finally, it is emphasized again that our goal was not to propose a method in the sense of introducing a method in a competition for the fastest running time. Rather, its purpose was to propose an alternative procedure for finding a strictly complementary solution. Implementation issue of our algorithm and comparisons with other existing methods are future research subjects.

## Appendix A

**Proof of Theorem 2.1**

As an immediate consequence of Lemma 2.1, $C_X$ is a blunt cone. Hence, we prove only that the convexity of $X$ is a necessary and sufficient condition for the convexity of $C_X$.

(*Sufficiency*) Assume that $X$ is convex and let $\mathbf{z}', \mathbf{z}'' \in C_X$. Then, by part (i) of Lemma 2.1, there exists $\alpha', \alpha'' > 0$ and $\mathbf{x}', \mathbf{x}'' \in X$ such that $\mathbf{z}' = \alpha' \begin{pmatrix} \mathbf{x}' \\ 1 \end{pmatrix}$ and $\mathbf{z}'' = \alpha'' \begin{pmatrix} \mathbf{x}'' \\ 1 \end{pmatrix}$. For any $\lambda \in (0,1)$, we define

$$\mathbf{z}_\lambda := \lambda \mathbf{z}' + (1-\lambda) \mathbf{z}'', \ \theta_\lambda := \lambda \frac{\alpha'}{\alpha''} + (1-\lambda), \ \hat{\alpha} := \alpha'' \theta_\lambda. \tag{A.1}$$

Then, $\mathbf{z}_\lambda$ can be expressed equivalently as

$$\mathbf{z}_\lambda = \begin{pmatrix} \lambda \alpha' \mathbf{x}' + (1-\lambda) \alpha'' \mathbf{x}'' \\ \lambda \alpha' + (1-\lambda) \alpha'' \end{pmatrix} = \begin{pmatrix} \lambda \frac{\alpha'}{\alpha''} (\hat{\alpha} \mathbf{x}') + \frac{1-\lambda}{\theta_\lambda} (\hat{\alpha} \mathbf{x}'') \\ \hat{\alpha} \end{pmatrix}. \tag{A.2}$$

By part (ii) of Lemma 2.1, the convexity of $X$ implies that $X^{\hat{\alpha}}$ is convex and so $\lambda \frac{\alpha'}{\alpha''} (\hat{\alpha} \mathbf{x}') + \frac{1-\lambda}{\theta_\lambda} (\hat{\alpha} \mathbf{x}'') \in X^{\hat{\alpha}}$. Hence, there exists $\hat{\mathbf{x}} \in X$ such that $\mathbf{z}_\lambda = \hat{\alpha} \begin{pmatrix} \hat{\mathbf{x}} \\ 1 \end{pmatrix}$ which indicates that $\mathbf{z}_\lambda \in C_X$.

(*Necessity*) Assume that $C_X$ is convex and let $\mathbf{x}_1, \mathbf{x}_2 \in X$. Then, $\begin{pmatrix} \mathbf{x}_1 \\ 1 \end{pmatrix}$ and $\begin{pmatrix} \mathbf{x}_2 \\ 1 \end{pmatrix}$ belong to $C_X$ by (2). Hence, by the assumption, $\begin{pmatrix} \lambda \mathbf{x}_1 + (1-\lambda) \mathbf{x}_2 \\ 1 \end{pmatrix} \in C_X$ for any $\lambda \in (0,1)$. Again by (2), it is followed that $X$ is convex. ∎

**Proof of Theorem 2.2**

(i) Let $\alpha > 0$ be given. Further, assume that $\hat{\mathbf{d}} \neq \mathbf{0}_n$ is an extreme direction of $X$ and that $\alpha \begin{pmatrix} \mathbf{x} \\ 1 \end{pmatrix} \in X^\alpha$. Then, by Definition 2.4, $\mathbf{x} + \frac{\theta}{\alpha} \hat{\mathbf{d}} \in X$ for all $\theta \geq 0$. Hence, $\alpha \begin{pmatrix} \mathbf{x} \\ 1 \end{pmatrix} + \theta \begin{pmatrix} \hat{\mathbf{d}} \\ 0 \end{pmatrix} = \alpha \begin{pmatrix} \mathbf{x} + \frac{\theta}{\alpha} \hat{\mathbf{d}} \\ 1 \end{pmatrix} \in X^\alpha$, showing that $\begin{pmatrix} \hat{\mathbf{d}} \\ 0 \end{pmatrix}$ is a recession direction of $X^\alpha$.

Conversely, let $\begin{pmatrix} \hat{\mathbf{d}} \\ 0 \end{pmatrix} \neq \mathbf{0}_{n+1}$ be a recession direction of $X^\alpha$ for all $\alpha > 0$ and let $\mathbf{x} \in X$. Then, $\hat{\mathbf{d}} \neq \mathbf{0}_n$ and $\begin{pmatrix} \mathbf{x} \\ 1 \end{pmatrix} \in X^1$ by (2). In addition, by definition 2.4, $\begin{pmatrix} \mathbf{x} + \theta \hat{\mathbf{d}} \\ 1 \end{pmatrix} \in X^1$ for all $\theta \geq 0$. Hence, $\mathbf{x} + \theta \hat{\mathbf{d}} \in X$ by (2) which indicates that $\hat{\mathbf{d}}$ is a recession direction of $X$.

Now, we prove that $\begin{pmatrix} \hat{\mathbf{d}} \\ 0 \end{pmatrix}$ is extreme if and only if $\hat{\mathbf{d}}$ is extreme. By Definition 2.4, $\begin{pmatrix} \hat{\mathbf{d}} \\ 0 \end{pmatrix}$ is an extreme direction if and only if it cannot be expressed as a positive combination of two distinct recession directions $\begin{pmatrix} \mathbf{d}_1 \\ 0 \end{pmatrix}$ and $\begin{pmatrix} \mathbf{d}_2 \\ 0 \end{pmatrix}$. This holds if and only if $\hat{\mathbf{d}}$ cannot be expressed as a positive combination of the two distinct recession directions $\mathbf{d}_1$ and $\mathbf{d}_2$; or, equivalently, if and only if $\hat{\mathbf{d}}$ is extreme.

**(ii)** By Remark 2.1, the convexity of $F$ is necessary and sufficient for the convexity of $F^\alpha$ for all $\alpha > 0$. Hence, let $F$ is a face of $X$ and let $\alpha \begin{pmatrix} \mathbf{x}' \\ 1 \end{pmatrix}, \alpha \begin{pmatrix} \mathbf{x}'' \\ 1 \end{pmatrix} \in X^\alpha$ with $\alpha > 0$ and $\mathbf{x}', \mathbf{x}'' \in X$. We define $\mathbf{z}_\lambda := \lambda \alpha \begin{pmatrix} \mathbf{x}' \\ 1 \end{pmatrix} + (1-\lambda) \alpha \begin{pmatrix} \mathbf{x}'' \\ 1 \end{pmatrix}$ for $\lambda \in (0,1)$. Then, $\mathbf{z}_\lambda$ can be expressed equivalently as $\mathbf{z}_\lambda = \alpha \begin{pmatrix} \lambda \mathbf{x}' + (1-\lambda) \mathbf{x}'' \\ 1 \end{pmatrix}$. By Definition 2.1, $F^\alpha$ is a face of $X^\alpha$ if $\mathbf{z}_\lambda \in F^\alpha$ with $\lambda \in (0,1)$ implies that $\alpha \begin{pmatrix} \mathbf{x}' \\ 1 \end{pmatrix}, \alpha \begin{pmatrix} \mathbf{x}'' \\ 1 \end{pmatrix} \in F^\alpha$. Assuming that $\mathbf{z}_\lambda \in F^\alpha$, we have $\lambda \mathbf{x}' + (1-\lambda) \mathbf{x}'' \in F$ and, by the assumption, we conclude that $\mathbf{x}', \mathbf{x}'' \in F$. Therefore, $\alpha' \begin{pmatrix} \mathbf{x}' \\ 1 \end{pmatrix}, \alpha'' \begin{pmatrix} \mathbf{x}'' \\ 1 \end{pmatrix} \in F^\alpha$ which indicates that $F^\alpha$ is a face of $X^\alpha$.

Conversely, assume that $F^\alpha$ is a face of $X^\alpha$ for all $\alpha > 0$. For $\lambda \in (0,1)$, we define $\mathbf{x}_\lambda := \lambda \mathbf{x}' + (1-\lambda) \mathbf{x}''$ with $\mathbf{x}', \mathbf{x}'' \in X$, and let $\mathbf{x}_\lambda \in F$ for some $\lambda \in (0,1)$. Then, by Definition 2.1, $F$ will be a face of $X$ if we show that $\mathbf{x}', \mathbf{x}'' \in F$. We have $\begin{pmatrix} \mathbf{x}_\lambda \\ 1 \end{pmatrix} \in F^1$, $\begin{pmatrix} \mathbf{x}' \\ 1 \end{pmatrix}, \begin{pmatrix} \mathbf{x}'' \\ 1 \end{pmatrix} \in X^1$ and $\begin{pmatrix} \mathbf{x}_\lambda \\ 1 \end{pmatrix} = \lambda \begin{pmatrix} \mathbf{x}' \\ 1 \end{pmatrix} + (1-\lambda) \begin{pmatrix} \mathbf{x}'' \\ 1 \end{pmatrix}$. Then, the assumption follows that $\begin{pmatrix} \mathbf{x}' \\ 1 \end{pmatrix}, \begin{pmatrix} \mathbf{x}'' \\ 1 \end{pmatrix} \in F^1$ and so $\mathbf{x}', \mathbf{x}'' \in F$.

∎

**Proof of Theorem 2.3**

By Remark 2.1, the convexity of $F$ is a necessary and sufficient condition for the convexity of $C_F$. So, let $F$ is a face of $X$ and let $\alpha' \begin{pmatrix} \mathbf{x}' \\ 1 \end{pmatrix}, \alpha'' \begin{pmatrix} \mathbf{x}'' \\ 1 \end{pmatrix} \in C_X$. For any $\lambda \in (0,1)$, we define $\mathbf{z}_\lambda := \lambda \alpha' \begin{pmatrix} \mathbf{x}' \\ 1 \end{pmatrix} + (1-\lambda) \alpha'' \begin{pmatrix} \mathbf{x}'' \\ 1 \end{pmatrix}$. Then, $C_F$ is a face of $C_X$ if $\mathbf{z}_\lambda \in C_F$ with $\lambda \in (0,1)$ implies that $\alpha' \begin{pmatrix} \mathbf{x}' \\ 1 \end{pmatrix}, \alpha'' \begin{pmatrix} \mathbf{x}'' \\ 1 \end{pmatrix} \in C_F$. By defining $\theta_\lambda$ and $\hat{\alpha}$ as those in (A.1), we express $\mathbf{z}_\lambda$ equivalently as that in (A.2) which implies the existence of a vector $\hat{\mathbf{x}} \in X$ such that $\mathbf{z}_\lambda = \hat{\alpha} \begin{pmatrix} \hat{\mathbf{x}} \\ 1 \end{pmatrix}$. Moreover, $\mathbf{z}_\lambda \in F^{\hat{\alpha}}$ since $\mathbf{z}_\lambda \in C_F$ and $F^{\hat{\alpha}} = X^{\hat{\alpha}} \cap C_F$. By the assumption and part (ii) of Theorem 2.2, $F^{\hat{\alpha}}$ is a face of $X^{\hat{\alpha}}$. Therefore, we conclude that $\alpha' \begin{pmatrix} \mathbf{x}' \\ 1 \end{pmatrix}, \alpha'' \begin{pmatrix} \mathbf{x}'' \\ 1 \end{pmatrix} \in F^{\hat{\alpha}}$, which indicates that $C_F$ is a face of $C_X$.

Conversely, assume that $C_F$ is a face of $C_X$. For $\lambda \in (0,1)$, let us define $\mathbf{x}_\lambda := \lambda \mathbf{x}' + (1-\lambda) \mathbf{x}''$ with $\mathbf{x}', \mathbf{x}'' \in X$. Further, let $\mathbf{x}_\lambda \in F$ for some $\lambda \in (0,1)$. Then, $F$ is a face of $X$, if we show that $\mathbf{x}', \mathbf{x}'' \in F$. Since $\begin{pmatrix} \mathbf{x}_\lambda \\ 1 \end{pmatrix} \in F^1$, by the assumption and part (ii) of Theorem 2.2, we have $\begin{pmatrix} \mathbf{x}' \\ 1 \end{pmatrix}, \begin{pmatrix} \mathbf{x}'' \\ 1 \end{pmatrix} \in F^1$ and so $\mathbf{x}', \mathbf{x}'' \in F$.

Now, it remains to prove the dimensionality statement. Let $F$ be a face of $X$ of dimension $k$. Then, there exists $k+1$ affinely independent $\mathbf{x}^1, ..., \mathbf{x}^{k+1}$ in $\text{aff}(F)$. Obviously, the vectors $\mathbf{0}_{n+1}, \begin{pmatrix} \mathbf{x}^1 \\ 1 \end{pmatrix}, ..., \begin{pmatrix} \mathbf{x}^{k+1} \\ 1 \end{pmatrix}$ are in $\text{aff}(C_F)$. The dimension of $C_F$ equals to $k+1$ if we show that

$\begin{pmatrix} \mathbf{x}^1 \\ 1 \end{pmatrix}, ..., \begin{pmatrix} \mathbf{x}^{k+1} \\ 1 \end{pmatrix}$ are linearly independent. For this purpose, let $\sum_{j=1}^{k+1} \alpha_j \begin{pmatrix} \mathbf{x}^j \\ 1 \end{pmatrix} = 0$ with $\boldsymbol{\alpha} \in \mathbb{R}^{k+1}$.

Then, $\sum_{j=2}^{k+1} \alpha_j (\mathbf{x}^j - \mathbf{x}^1) = 0$ and $\sum_{j=1}^{k+1} \alpha_j = 0$. Since $\mathbf{x}^2 - \mathbf{x}^1, ..., \mathbf{x}^{k+1} - \mathbf{x}^1$ are linearly independency, we have $\boldsymbol{\alpha} = \mathbf{0}_{k+1}$ which indicates that $C_F$ is of dimension $k+1$.

Conversely, let $C_F$ is of dimension $k+1$. Then, there exists $k+2$ affinely independent $\begin{pmatrix} \mathbf{x}^1 \\ 1 \end{pmatrix}, ..., \begin{pmatrix} \mathbf{x}^{k+2} \\ 1 \end{pmatrix}$ in aff$(C_F)$. This follows that the vectors $\mathbf{x}^2 - \mathbf{x}^1, ..., \mathbf{x}^{k+2} - \mathbf{x}^1$ in aff$(F)$ are linearly independent and so $F$ is of dimension $k$. ∎

**Proof of Theorem 2.4**

Let $\begin{pmatrix} \hat{\mathbf{d}} \\ \hat{d}_{n+1} \end{pmatrix}$ be a recession direction of $C_X$ and let $\mathbf{x} \in X$ and $\theta \geq 0$ be arbitrary. Then, by Definition 2.4, $\alpha \begin{pmatrix} \mathbf{x} \\ 1 \end{pmatrix} + \theta \begin{pmatrix} \hat{\mathbf{d}} \\ \hat{d}_{n+1} \end{pmatrix} \in C_X$ for all $\alpha > 0$. If $\hat{d}_{n+1} = 0$, by setting $\alpha = 1$, we obtain $\mathbf{x} + \theta \hat{\mathbf{d}} \in X$ which implies that $\hat{\mathbf{d}}$ is a recession direction of $X$. If $\hat{d}_{n+1} > 0$, we define $\lambda := \frac{\alpha}{\alpha + \theta \hat{d}_{n+1}}$. Then, it is straightforward that $\lambda \mathbf{x} + (1 - \lambda) \frac{\hat{\mathbf{d}}}{\hat{d}_{n+1}} \in X$. Obviously, if $\alpha \to 0$ then $\lambda \to 0$ and consequently $\frac{\hat{\mathbf{d}}}{\hat{d}_{n+1}} \in X$ by the closedness of $X$.

Conversely, let $\hat{\mathbf{d}}$ be a recession direction of $X$ and let $\alpha > 0$, $\theta \geq 0$ and $\mathbf{x} \in X$ be arbitrary. By definition, we have $\mathbf{x} + \frac{\theta}{\alpha} \hat{\mathbf{d}} \in X$ and so $\alpha \begin{pmatrix} \mathbf{x} \\ 1 \end{pmatrix} + \theta \begin{pmatrix} \hat{\mathbf{d}} \\ 0 \end{pmatrix} \in C_X$. This implies that $\begin{pmatrix} \hat{\mathbf{d}} \\ 0 \end{pmatrix}$ is a recession direction of $C_X$.

Since any extreme points and any extreme direction of a convex set are its 0- and 1-faces, the proof is complete by Theorem 2.4. ∎

**Proof of Theorem 2.5**

**(i)** ($\Rightarrow$) Let $X$ be unbounded. Then, $D_h \neq \emptyset$. First, we show that $C_X \cup coni(D_h) \subseteq cl(C_X)$. Since $cl(C_X)$ is a convex cone containing $C_X$, it is enough to show only that $D_h \subseteq cl(C_X)$. Let $\begin{pmatrix} \hat{\mathbf{d}} \\ 0 \end{pmatrix} \in D_h$ with $\|\hat{\mathbf{d}}\| = 1$ and let $\mathbf{x} \in X$. Then, $\frac{1}{n}\begin{pmatrix} \mathbf{x} + n\hat{\mathbf{d}} \\ 1 \end{pmatrix}$ is a sequence in $C_X$ whose convergence to $\begin{pmatrix} \hat{\mathbf{d}} \\ 0 \end{pmatrix}$ follows that $\begin{pmatrix} \hat{\mathbf{d}} \\ 0 \end{pmatrix} \in cl(C_X)$.

To prove the inverse inclusion, let $\hat{\mathbf{z}} \in cl(C_X)$. If $\hat{\mathbf{z}} = \mathbf{0}_{n+1}$, it is obvious that $\hat{\mathbf{z}} \in coni(D_h)$ and so the result is obtained. If $\hat{\mathbf{z}} \neq \mathbf{0}_{n+1}$, we can construct the sequence $\{\mathbf{z}_j\}$ in $C_X$ so that it converges to $\hat{\mathbf{z}}$. For any $j \in \mathbb{N}$, there exist $\mathbf{x}_j \in X$ and $\alpha_j > 0$ such that $\mathbf{z}_j = \alpha_j \begin{pmatrix} \mathbf{x}_j \\ 1 \end{pmatrix}$. Hence, $\alpha_j \mathbf{x}_j$ and $\alpha_j$ converge to $(\hat{z}_1, ..., \hat{z}_n)^T$ and $\hat{z}_{n+1}$, respectively. We now consider two cases:

*Case 1:* $\hat{z}_{n+1} > 0$. Then, the sequence $\{\mathbf{x}_j\}$ converges to $\frac{1}{\hat{z}_{n+1}}(\hat{z}_1, ..., \hat{z}_n)^T$. Since $X$ is closed, we have $\frac{1}{\hat{z}_{n+1}}(\hat{z}_1, ..., \hat{z}_n)^T \in X$ thereby $\hat{\mathbf{z}} \in X^{\hat{z}_{n+1}} \subseteq C_X$.

*Case 2:* $\hat{z}_{n+1} = 0$. Then, for any arbitrary $\mathbf{x} \in X$ and $\theta \geq 0$, the convexity of $X$ implies that the sequence $\{\bar{\mathbf{x}}_j\}$ defined by $\bar{\mathbf{x}}_j := \frac{1}{1+\theta\alpha_j}\mathbf{x} + \frac{\theta\alpha_j}{1+\theta\alpha_j}\mathbf{x}_j$, $j \in \mathbb{N}$, lies in $X$. Since $X$ is closed, the limit of this sequence belongs to $X$, i.e., $\mathbf{x} + \theta(\hat{z}_1, ..., \hat{z}_n)^T \in X$. This means that $(\hat{z}_1, ..., \hat{z}_n)^T$ is a recession direction of $X$ and so $\hat{\mathbf{z}} \in coni(D_h)$.

($\Leftarrow$) Let $cl(C_X) = C_X \cup coni(D_h)$. We claim that $X$ is unbounded. Assume on the contrary that $X$ is bounded. Then, $D_h = \emptyset$ which implies that $cl(C_X) = C_X$ which is a contradiction.

**(ii)** Let $X$ be bounded. Since $C_X \cup \{\mathbf{0}_{n+1}\} \subseteq cl(C_X)$, we prove only the inverse inclusion. For this purpose, let $\hat{\mathbf{z}} \in cl(C_X)$. Then, there exists a sequence $\{\mathbf{z}_j\}$ in $C_X$ which converges to $\hat{\mathbf{z}}$. For any $j \in \mathbb{N}$, there exist $\mathbf{x}_j \in X$ and $\alpha_j \geq 0$ such that $\mathbf{z}_j = \alpha_j \begin{pmatrix} \mathbf{x}_j \\ 1 \end{pmatrix}$. The boundedness of $X$ then follows that the sequence $\{\mathbf{x}_j\} \subseteq X$ has a convergent subsequence, namely $\{\mathbf{x}_{j_k}\}$. Assuming that $\hat{\mathbf{x}}$ be the limit of this subsequence, the closedness of $X$ implies that $\hat{\mathbf{x}} \in X$. Therefore, the subsequence $\{\mathbf{z}_{j_k}\}$ converges to $\hat{z}_{n+1}\begin{pmatrix}\hat{\mathbf{x}}\\1\end{pmatrix}$. Since the limit of a sequence is the same that of each of its subsequences, we have $\hat{\mathbf{z}} = \hat{z}_{n+1}\begin{pmatrix}\hat{\mathbf{x}}\\1\end{pmatrix}$ indicating that $cl(C_X) \subseteq C_X \cup \{\mathbf{0}_{n+1}\}$.

Conversely, if $cl(C_X) = C_X \cup \{\mathbf{0}_{n+1}\}$, then part (i) of the theorem implies that $X$ is bounded.

∎

**Proof of Theorem 2.6**

By Theorem 2.4 and (3), it is enough to show only that $C_X$ and $cl(C_X)$ have the same recession directions. Let $\begin{pmatrix}\hat{\mathbf{d}}\\\hat{d}_{n+1}\end{pmatrix}$ be a recession direction of $cl(C_X)$ and let $\theta \geq 0$ be arbitrary. By definition, we have $\mathbf{z} + \theta\begin{pmatrix}\hat{\mathbf{d}}\\\hat{d}_{n+1}\end{pmatrix} \in cl(C_X)$ for any $\mathbf{z} \in cl(C_X)$. To prove that $\begin{pmatrix}\hat{\mathbf{d}}\\\hat{d}_{n+1}\end{pmatrix}$ is a recession direction of $C_X$, let $\mathbf{z} \in C_X$. Then, Theorem 2.5 implies that $\mathbf{z} + \theta\begin{pmatrix}\hat{\mathbf{d}}\\\hat{d}_{n+1}\end{pmatrix} \in C_X$ since $z_{n+1} + \theta\hat{d}_{n+1} > 0$. This indicates that $\begin{pmatrix}\hat{\mathbf{d}}\\\hat{d}_{n+1}\end{pmatrix}$ is a recession direction of $C_X$.

Conversely, let $\begin{pmatrix}\hat{\mathbf{d}}\\\hat{d}_{n+1}\end{pmatrix}$ be a recession direction of $C_X$ and let $\theta \geq 0$ be arbitrary. By definition, we have $\mathbf{z} + \theta\begin{pmatrix}\hat{\mathbf{d}}\\\hat{d}_{n+1}\end{pmatrix} \in C_X$ for any $\mathbf{z} \in C_X$. To prove that $\begin{pmatrix}\hat{\mathbf{d}}\\\hat{d}_{n+1}\end{pmatrix}$ is a recession direction of $cl(C_X)$, let $\mathbf{z} \in cl(C_X)$ and consider two cases:

*Case 1:* $X$ is unbounded. Then, $D_h \neq \emptyset$ and part (i) of Theorem 2.5 follows that $\mathbf{z}$ belongs to either $C_X$ or $coni(D_h)$.

If $\mathbf{z} \in C_X$, the result is direct. So, let $\mathbf{z} \in coni(D_h)$. If $\hat{d}_{n+1} = 0$, then part (i) of Theorem 2.4 follows that $\begin{pmatrix} \hat{\mathbf{d}} \\ \hat{d}_{n+1} \end{pmatrix} \in coni(D_h)$ and so the result is immediate by the convexity of $coni(D_h)$. So, let $\hat{d}_{n+1} > 0$. Then, part (ii) of Theorem 2.4 and the convexity of $cl(C_X)$ implies that $\mathbf{z} + \theta \begin{pmatrix} \hat{\mathbf{d}} \\ \hat{d}_{n+1} \end{pmatrix} \in cl(C_X)$. This indicates that $\begin{pmatrix} \hat{\mathbf{d}} \\ \hat{d}_{n+1} \end{pmatrix}$ is a recession direction of $cl(C_X)$.

*Case 2:* $X$ is bounded. Then, $D_h = \emptyset$ and so $d_{n+1} > 0$. In addition, part (ii) of Theorem 2.5 follows that $\mathbf{z} \in C_X$ or $\mathbf{z} = \mathbf{0}_{n+1}$. If $\mathbf{z} \in C_X$, the result is direct. If $\mathbf{z} = \mathbf{0}_{n+1}$, the result is obtained by part (ii) of Theorem 2.4. Thus, the proof is completes. ∎

**Proof of Theorem 2.7**

By Theorem 2.6, $cone(conv(D_v))$ and $coni(D_h)$ are subsets of $cl(C_X)$ whose direct sum is also a subset of $cl(C_X)$ by the convexity of this cone. Hence, we prove only the inverse of inclusion. Based on the fact that a closed convex cone containing no line is the convex hull of its extreme directions (Fenchel, 1951), Theorem 2.6 implies that

$$cl(C_X) = conv(cone(D_v \cup D_h)). \tag{A.3}$$

Based on (A.3) and the Caratheodory theorem (Rockafellar, 1970), $\mathbf{z} \in cl(C_X)$ if and only if there exist $n+2$ vectors in $cone(D_v \cup D_h)$ such that $\mathbf{z}$ can be expressed as a convex combination of them. Formally, there exist $\alpha_1 \begin{pmatrix} \mathbf{d}^1 \\ 1 \end{pmatrix}, \dots, \alpha_k \begin{pmatrix} \mathbf{d}^k \\ 1 \end{pmatrix}, \alpha_{k+1} \begin{pmatrix} \mathbf{d}^{k+1} \\ 0 \end{pmatrix}, \dots, \alpha_{n+2} \begin{pmatrix} \mathbf{d}^{n+2} \\ 0 \end{pmatrix}$ $\in cone(D_v \cup D_h)$ with $\alpha_j \geq 0$, $j=1,\dots,n+2$, $\lambda_j \geq 0$, $j=1,\dots,n+2$, and $\sum_{j=1}^{n+2} \lambda_j = 1$ such that

$$\mathbf{z} = \sum_{j=1}^{k} \lambda_j \alpha_j \begin{pmatrix} \mathbf{d}^j \\ 1 \end{pmatrix} + \sum_{j=k+1}^{n+2} \lambda_j \alpha_j \begin{pmatrix} \mathbf{d}^j \\ 0 \end{pmatrix}. \tag{A.4}$$

If $z_{n+1} = 0$, then $\sum_{j=1}^{k} \lambda_j \alpha_j = 0$ implying that $\lambda_j \alpha_j = 0$ for all $j = 1,...,k$. In this case, $\sum_{j=1}^{k} \lambda_j \alpha_j \begin{pmatrix} \mathbf{x}^j \\ 1 \end{pmatrix} = \mathbf{0}_{n+1} \in cone(conv(D_v))$. Moreover, because of the convexity of $coni(D_h)$, we have $\sum_{j=k+1}^{n+2} \lambda_j \alpha_j \begin{pmatrix} \mathbf{d}^j \\ 0 \end{pmatrix} \in coni(D_h)$. Thus, $\mathbf{z} \in cone(conv(D_v)) + coni(D_h)$.

If $z_{n+1} > 0$, by defining $\gamma_j := \dfrac{\lambda_j \alpha_j}{z_{n+1}}$, $j = 1,...,n+2$, we then rewrite $\mathbf{z}$ as

$$\mathbf{z} = z_{n+1} \sum_{j=1}^{k} \gamma_j \begin{pmatrix} \mathbf{d}^j \\ 1 \end{pmatrix} + \sum_{j=k+1}^{n+2} z_{n+1} \gamma_j \begin{pmatrix} \mathbf{d}^j \\ 0 \end{pmatrix},$$
$$\sum_{j=1}^{k} \gamma_j = 1, \ \gamma_j \geq 0, \ j = 1,...,n+2,$$
(A.5)

which follows that $\mathbf{z} \in cone(conv(D_v)) + coni(D_h)$.

Hence, in both of the cases, we see that $cl(C_X) \subseteq cone(conv(D_v)) + coni(D_h)$ and so the proof is complete. ∎

**Proof of Lemma 3.1**

**(i)** This result is an immediate consequence of the convexity of $X$.

**(ii)** By Definition 3.3, we have $\Sigma(X) = \sigma(\bar{\mathbf{x}})$ where $\bar{\mathbf{x}}$ is a maximal element of $X$. Hence, we need only to prove that $\bigcup_{\mathbf{x} \in X} \sigma(\mathbf{x}) \subseteq \sigma(\bar{\mathbf{x}})$. That is, $\bar{\mathbf{x}}$ takes positive values in any positive components of any $\mathbf{x} \in X$. By way of contradiction, assume that there exists some vector $\hat{\mathbf{x}}$ in $X$ which takes positive values in some zero components of $\bar{\mathbf{x}}$. Without loss of generality, assume that $\sigma(\bar{\mathbf{x}}) = \{1,...,k\}$, and that $\{k+1,...,l\}$ is the index set of components in which $\hat{\mathbf{x}}$ takes positive values. Then, from Part (i), it follows that $\sigma(\bar{\mathbf{x}}) \subsetneq \sigma(\lambda \bar{\mathbf{x}} + (1-\lambda)\hat{\mathbf{x}})$ for any $\lambda \in (0,1)$, which contradicts $\bar{\mathbf{x}}$'s being maximal.

**(iii)** This result is an immediate consequence of part (ii) of the lemma.

**(iv)** The result is straightforward and so the proof is omitted. ∎

**Proof of Theorem 3.1**

First note that problem (6) is always feasible since $(\mathbf{z}^1, \mathbf{z}^2) = (\mathbf{0}_{n+1}, \mathbf{0}_{n+1})$ is a feasible solution for this problem. In addition, this problem has an optimal solution as the objective function is upper bounded by $n+1$. Hence, let $(\mathbf{z}^{1*}, \mathbf{z}^{2*})$ be an optimal solution to problem (6).

**(i)** Since $\mathbf{z}^{1*} \geq \mathbf{0}_{n+1}$, $z_j^{1*}$ takes either zero or positive value for every $j$. On the other hand, we have $\mathbf{z}^{1*} + \mathbf{z}^{2*} \geq \mathbf{0}_{n+1}$ since $\mathbf{z}^{1*} + \mathbf{z}^{2*} \in C$. Hence, $z_j^{2*} \geq 0$ for any $j$ that $z_j^{1*} = 0$. Moreover, since (6) is a maximization problem, $z_j^{1*}$ cannot take a positive value unless $z_j^{2*}$ be positive. Therefore, $z_j^{2*} > 0$ for any $j$ that $z_j^{1*} > 0$.

**(ii)** Assume by contradiction that there exists $l \in \{1, \ldots, n+1\}$ for which $0 < z_l^{2*} < 1$. Then, $z_l^{1*} = 0$ and the constraints of (6) at optimality can be expressed equivalently as

$$\left( z_1^{1*} + z_1^{2*}, \ldots, z_{l-1}^{1*} + z_{l-1}^{2*}, z_l^{2*}, z_{l+1}^{1*} + z_{l+1}^{2*}, \ldots, z_{n+1}^{1*} + z_{n+1}^{2*} \right)^T \in C. \tag{A.6}$$

Since $C$ is a cone and $z_l^{2*} > 0$, we have

$$\left( \frac{z_1^{1*} + z_1^{2*}}{z_l^{2*}}, \ldots, \frac{z_{l-1}^{1*} + z_{l-1}^{2*}}{z_l^{2*}}, 1, \frac{z_{l+1}^{1*} + z_{l+1}^{2*}}{z_l^{2*}}, \ldots, \frac{z_{n+1}^{1*} + z_{n+1}^{2*}}{z_l^{2*}} \right)^T \in C. \tag{A.7}$$

For $j = 1, \ldots, n+1$, we define the vector $(\mathbf{z}^{1\prime}, \mathbf{z}^{2\prime})$ as

$$z_j^{1\prime} := \max\left\{ 0, \frac{1}{z_l^{2*}}\left(z_j^{1*} + z_j^{2*}\right) - 1 \right\} \text{ and } z_j^{2\prime} := \min\left\{ 1, \frac{1}{z_l^{2*}}\left(z_j^{1*} + z_j^{2*}\right) \right\}. \tag{A.8}$$

Then, $z_j^{1\prime} = 0$ and $z_j^{2\prime} = \frac{1}{z_l^{2*}}\left(z_j^{1*} + z_j^{2*}\right)$ for any $j$ that $\frac{1}{z_l^{2*}}\left(z_j^{1*} + z_j^{2*}\right) \leq 1$. Moreover, $z_j^{1\prime} = \frac{1}{z_l^{2*}}\left(z_j^{1*} + z_j^{2*}\right) - 1$ and $z_j^{2\prime} = 1$ for any $j$ that $\frac{1}{z_l^{2*}}\left(z_j^{1*} + z_j^{2*}\right) > 1$. Thus, $z_j^{1\prime} + z_j^{2\prime} = \frac{1}{z_l^{2*}}\left(z_j^{1*} + z_j^{2*}\right)$ for every $j \neq l$ and $z_l^{1\prime} + z_l^{2\prime} = 1$. Therefore, according to (A.7), $(\mathbf{z}^{1\prime}, \mathbf{z}^{2\prime})$ is a feasible solution to (6) whose objective value is strictly greater than $\sum_{j=1}^{n+1} z_j^{2*}$. This contradicts the optimality of $(\mathbf{z}^{1*}, \mathbf{z}^{2*})$ and proves that $z_j^{2*} = 1$ for any $j$ that $z_j^{2*} > 0$.

**(iii)** The constraints of (6) guarantees that the vector $\bar{\mathbf{z}}$ defined by $\bar{\mathbf{z}} := \mathbf{z}^{1*} + \mathbf{z}^{2*}$ belongs to $C$. We claim that $\bar{\mathbf{z}}$ has the maximum number of positive components. Assume by contradiction that this is not true. Then, part (ii) of Lemma 3.1 follows that there exists a vector $\hat{\mathbf{z}} \leq \mathbf{1}_{n+1}$ in $C$ which takes positive values in some zero components of $\bar{\mathbf{z}}$. Without loss of generality, assume that $\sigma(\bar{\mathbf{z}}) = \{1,...,k\}$ and $\sum_{j=k+1}^{n+1} \hat{z}_j > 0$.

Since $\bar{\mathbf{z}}, \hat{\mathbf{z}} \in C$, the convexity of $C$ follows that $\bar{\mathbf{z}} + \hat{\mathbf{z}} \in C$ and consequently,

$$\left( \left( z_1^{1*} + \hat{z}_1 \right) + z_1^{2*}, ..., \left( z_k^{1*} + \hat{z}_k \right) + z_k^{2*}, \hat{z}_{k+1}, ..., \hat{z}_{n+1} \right)^T \in C. \tag{A.9}$$

Based on (A.9), we define

$$z_j^{1''} := \begin{cases} z_j^{1*} + \hat{z}_j, & j = 1,...,k, \\ 0, & j = k+1,...,n+1, \end{cases} \quad \text{and} \quad z_j^{2''} := \begin{cases} z_j^{2*}, & j = 1,...,k, \\ \hat{z}_j, & j = k+1,...,n+1. \end{cases} \tag{A.10}$$

Then, $\left( \mathbf{z}^{1''}, \mathbf{z}^{2''} \right)$ is a feasible solution to (6) with the objective value $\sum_{j=1}^{k} z_j^{2*} + \sum_{j=k+1}^{n+1} \hat{z}_j$, which is strictly greater than $\sum_{j=1}^{k} z_j^{2*}$. This contradicts the optimality of $\left( \mathbf{z}^{1*}, \mathbf{z}^{2*} \right)$ and proves our claim. ∎

**Proof of Theorem 3.2**

Let $\left( \begin{pmatrix} \mathbf{x}^{1*} \\ w^{1*} \end{pmatrix}, \begin{pmatrix} \mathbf{x}^{2*} \\ w^{2*} \end{pmatrix} \right)$ be an optimal solution to (7). First, we prove that $X = \emptyset$ if $w^{2*} = 0$.

Assume by contradiction that $\hat{\mathbf{x}} \in X$. Then, $\begin{pmatrix} \hat{\mathbf{x}} \\ 1 \end{pmatrix} \in C_X$. If we define $\mathbf{z}^{1'} := \begin{pmatrix} \hat{\mathbf{x}} \\ 0 \end{pmatrix}$ and $\mathbf{z}^{2'} := \begin{pmatrix} \mathbf{0}_n \\ 1 \end{pmatrix}$, then $\left( \mathbf{z}^{1'}, \mathbf{z}^{2'} \right)$ is a feasible solution to (7). Since the $(n+1)th$ component of $\mathbf{z}^{1'} + \mathbf{z}^{2'}$ is positive, part (ii) of Lemma 3.1 infers that $w^{2*} > 0$, which is a contradiction.

Now, let $w^{2*} > 0$. Then, parts (ii) and (iii) of Theorem 3.1 imply that $w^{2*} = 1$ and that $\begin{pmatrix} \mathbf{x}^{1*} + \mathbf{x}^{2*} \\ w^{1*} + 1 \end{pmatrix}$ is a maximal element of $C_X$. Since $C_X$ is a cone, $\begin{pmatrix} \bar{\mathbf{x}} \\ 1 \end{pmatrix}$ is also a maximal element of $C_X$. Thus, the proof is complete by part (iv) of Lemma 3.1. ∎

**Proof of Theorem 4.1**

3333

Let $\hat{\mathbf{x}} \in ri(P)$. From part (iii) of Lemma 3.1, we prove $\hat{\mathbf{x}} \in \arg\max_{\mathbf{x} \in P} \|\mathbf{x}\|_0$ by showing that $\hat{x}_j > 0$ for all $j \in J^{\neq} = \{j \mid x_j > 0 \text{ for some } \mathbf{x} \in X\}$. By contradiction, assume that $\hat{x}_{j_0} = 0$ for some $j_0 \in J^{\neq}$. In order to derive a contradiction, we demonstrate that $N_\varepsilon(\hat{\mathbf{x}}) \cap \text{aff}(P)$ is not contained in $X$ for any arbitrary $\varepsilon > 0$. From the definition of $J^{\neq}$, there exists $\overline{\mathbf{x}} \in P$ for which $\overline{x}_{j_0} > 0$. For $\delta > 0$, if we define $\mathbf{x}_\delta := (1+\delta)\hat{\mathbf{x}} - \delta\overline{\mathbf{x}}$, then $\mathbf{x}_\delta \in \text{aff}(P)$. Since $\delta$ is an arbitrary positive number, selecting $\delta_\varepsilon > 0$ as a sufficiently small number associated with $\varepsilon$, $\mathbf{x}_{\delta_\varepsilon}$ belongs to $N_\varepsilon(\hat{\mathbf{x}}) \cap \text{aff}(P)$. Since $\hat{x}_{j_0} = 0$ and $\overline{x}_{j_0} > 0$, we have $x_{\delta j_0} < 0$ indicating that $\mathbf{x}_\delta \notin P$. This is a contradiction and thus $ri(P) \subseteq \arg\max_{\mathbf{x} \in P} \|\mathbf{x}\|_0$.

To prove the reverse inclusion, let $\hat{\mathbf{x}} \in \arg\max_{\mathbf{x} \in P} \|\mathbf{x}\|_0$. From part (iii) of Lemma 3.1, we have $\hat{x}_j > 0$ for all $j \in J^{\neq}$. Then, it is easy to assert that there exists $\varepsilon > 0$ such that $x_j > 0$ for all $j \in J^{\neq}$ and for all $\mathbf{x}$ in $N_\varepsilon(\hat{\mathbf{x}})$. On the other hand, for all $\mathbf{x}$ in $\text{aff}(P)$, we have $\mathbf{A}\mathbf{x} = \mathbf{b}$ and $x_j = 0$ for all $j \in \{j \mid \hat{x}_j = 0\}$. Hence, $N_\varepsilon(\hat{\mathbf{x}}) \cap \text{aff}(P) \subseteq P$. This means that $\hat{\mathbf{x}} \in ri(P)$ and so the proof is complete. ∎

## Appendix B

The computer program written using GAMS:

```
1   Sets
2       i    row number of matrix A      /i1*i5/
3       j    column number of matrix A   /j1*j4/;
4
5   Table A(i,j)
6            j1       j2       j3       j4
7   i1       -1       1        -2       1
8   i2       4        -4       1        -2
9   i3       0        0        -3       1
10  i4       -1       1        -2       1
11  i5       -2       5        -9       3    ;
12
13  Parameters
14      b/i1    1
15       i2     0
16       i3     2
17       i4     1
18       i5     7 /
19      c/j1    -4
20       j2     4
21       j3     -8
22       j4     4 /
23      Zstar
24      Xstar(j)
25      Ustar(i)
```

```
26        Ystar(i)
27        Vstar(j);
28
29  Free Variables
30        Teta;
31
32  Positive Variables
33        x(j),t(j),u(i),s(i),x1,t1,y(i),p(i),v(j),q(j),y1,p1;
34
35  t.up(j) = 1;
36  s.up(i) = 1;
37  p.up(i) = 1;
38  q.up(j) = 1;
39  t1.up = 1;
40  p1.up = 1;
41
42  Equations
43        Obj
44          Con
45        ObjP
46          ConP1
47          ConP2
48        ObjD
49          ConD1
50          ConD2
51        Obj1
52          ConP
53          ConD
54          ConOpt ;
55
56        Obj..         Teta =E= Sum(j, c(j)*x(j));
57         Con(i)..        Sum(j, a(i,j)*x(j)) =L= b(i);
58
59        ObjP..        Teta =E= Sum(j, t(j)) + t1 + Sum(i, s(i));
60          ConP1(i)..      Sum(j, a(i,j)*(x(j)+t(j)) ) + u(i)+s(i) -  b(i)*(x1+t1) =E= 0;
61          ConP2..         Sum(j, c(j)  *(x(j)+t(j)) )                - Zstar*(x1+t1) =E= 0;
62
63        ObjD..        Teta =E= Sum(i, p(i)) + p1 + Sum(j, q(j));
64          ConD1(j)..      Sum(i, a(i,j)*(y(i)+p(i)) ) - v(j)-q(j) -  c(j)*(y1+p1) =E= 0;
65          ConD2..         Sum(i, b(i)  *(y(i)+p(i)) )                - Zstar*(y1+p1) =E= 0;
66
67        Obj1..    Teta =E= Sum(j, t(j)) + t1 + Sum(i, s(i)) + Sum(i, p(i)) + Sum(j,q(j));
68          ConP(i)..       Sum(j, a(i,j)*(x(j)+t(j)) ) + u(i)+s(i) -  b(i)*(x1+t1) =E= 0;
69          ConD(j)..       Sum(i, a(i,j)*(y(i)+p(i)) ) - v(j)-q(j) -  c(j)*(x1+t1) =E= 0;
70          ConOpt..        Sum(j, c(j)*(x(j)+t(j)) ) - Sum(i, b(i)*(y(i)+p(i)) )   =E= 0;
71
72  Model MainLP       / Obj, Con/
73        Primal_SCSC  / ObjP, ConP1, ConP2/
74        Dual_SCSC    / ObjD, ConD1, ConD2/
75        SCSC         / Obj1, ConP, ConD, ConOpt/;
76
77
78  Solve MainLP using LP Maximizing Teta;
79        Zstar=Teta.L;
80
81  Solve Primal_SCSC using LP Maximizing Teta;
82  Solve   Dual_SCSC using LP Maximizing Teta;
83        Xstar(j) = (x.L(j)+t.L(j))/(x1.L+t1.L);
84        Ustar(i) = (u.L(i)+s.L(i))/(x1.L+t1.L);
85        Ystar(i) = (y.L(i)+p.L(i))/(y1.L+p1.L);
86        Vstar(j) = (v.L(j)+q.L(j))/(y1.L+p1.L);
87
88  Display Xstar, Vstar, Ustar, Ystar;
89
90  Solve SCSC using LP Maximizing Teta;
91        Xstar(j) = (x.L(j)+t.L(j))/(x1.L+t1.L);
92        Ustar(i) = (u.L(i)+s.L(i))/(x1.L+t1.L);
93        Ystar(i) = (y.L(i)+p.L(i))/(x1.L+t1.L);
94        Vstar(j) = (v.L(j)+q.L(j))/(x1.L+t1.L);
95
96  Display Xstar, Vstar, Ustar, Ystar;
97
98
```